\documentclass{amsart}

\usepackage{graphicx}       % standard LaTeX graphics tool
\usepackage[latin1]{inputenc}
\usepackage[T1]{fontenc}
\usepackage{microtype} % good font tricks

\usepackage{array,colortbl}
\usepackage{amsmath,amsfonts,amssymb,bm,amsthm,mathtools,mathrsfs}
\DeclareFontFamily{U}{mathx}{\hyphenchar\font45}
\DeclareFontShape{U}{mathx}{m}{n}{
      <5> <6> <7> <8> <9> <10>
      <10.95> <12> <14.4> <17.28> <20.74> <24.88>
      mathx10
      }{}
\DeclareSymbolFont{mathx}{U}{mathx}{m}{n}
\DeclareFontSubstitution{U}{mathx}{m}{n}
\DeclareMathAccent{\widecheck}{0}{mathx}{"71}

\def\citep#1#2{\cite[{#1}]{#2}}

\newcommand{\RefEq}[1]{~\textup{(\ref{#1})}}
\newcommand{\RefEqTwo}[2]{~\textup{(\ref{#1})} and~\textup{(\ref{#2})}}

\newcommand{\RefSec}[1]{Section~\textup{\ref{#1}}}

\newcommand{\RefThm}[1]{Theorem~\textup{\ref{#1}}}

\newcommand{\RefLem}[1]{Lemma~\textup{\ref{#1}}}

\newcommand{\RefCol}[1]{Corollary~\textup{\ref{#1}}}

\newcommand{\RefAlg}[1]{Algorithm~\textup{\ref{#1}}}

\newcommand{\RefFig}[1]{Figure~\textup{\ref{#1}}}

\usepackage[ruled,linesnumbered]{algorithm2e}

\newcommand{\N}{{\mathbb{N}}} % natural numbers {1, 2, ...}
\newcommand{\R}{{\mathbb{R}}} % reals
\newcommand{\Z}{{\mathbb{Z}}} % integers
\newcommand{\bsh}{{\boldsymbol{h}}}

\newcommand{\bsk}{{\boldsymbol{k}}}

\newcommand{\bss}{{\boldsymbol{s}}}
\newcommand{\bst}{{\boldsymbol{t}}}
\newcommand{\bsu}{{\boldsymbol{u}}}
\newcommand{\bsv}{{\boldsymbol{v}}}

\newcommand{\bsx}{{\boldsymbol{x}}}

\newcommand{\bsz}{{\boldsymbol{z}}}

\newcommand{\bsell}{{\boldsymbol{\ell}}}
\newcommand{\bszero}{{\boldsymbol{0}}} % vector of zeros
\newcommand{\bsone}{{\boldsymbol{1}}}  % vector of ones
\newcommand{\bsgamma}{{\boldsymbol{\gamma}}}

\newcommand{\bstheta}{{\boldsymbol{\theta}}}

\newcommand{\calA}{{\mathcal{A}}}

\newcommand{\calC}{{\mathcal{C}}}

\newcommand{\calX}{{\mathcal{X}}}
\newcommand{\calY}{{\mathcal{Y}}}
\newcommand{\calZ}{{\mathcal{Z}}}

\newcommand{\frakg}{{\mathfrak{g}}}
\newcommand{\frakh}{{\mathfrak{h}}}

\newcommand{\abs}[1]{\left\vert#1\right\vert}
\newcommand{\norm}[1]{\left\Vert#1\right\Vert}
\newcommand{\floor}[1]{\left\lfloor #1 \right\rfloor} % floor
\newcommand{\ceil}[1]{\left\lceil #1 \right\rceil}    % ceil
\newcommand{\rd}{\,\mathrm{d}} % differential symbol with tiny space in front for use in integrals
\newcommand{\innProd}[3]{\left\langle #1, #2 \right\rangle_{#3}}

\newcommand{\gmu}{{\gamma_{\bsu}}}

\newcommand{\fh}{\widehat{f}}

\newcommand{\WXZl}{W_{\calX,\calZ,\nu,\ell}}
\newcommand{\WXYZl}{W_{\calX,\calY,\calZ,\nu,\ell}}

\newcommand{\tauup}{\tau^{\textup{up}}}
\newcommand{\taulow}{\tau^{\textup{low}}}

\newcommand{\ghk}{\widehat{g}_{\bsk}}

\newcommand{\Halphad}{H_{\alpha,d}}

\newcommand{\convRate}{\alpha - \frac{1}{2} - \delta}

\newcommand{\Halphaminusgamma}{H_{\convRate,\bsgamma,d}}
\newcommand{\Knugammad}{K_{\nu,\bsgamma,d}^{\alpha}}
\newcommand{\Qmgammad}{Q_{m,\bsgamma,d}^{\alpha}}
\newcommand{\Knugammau}{K_{\nu,\bsgamma,\bsu}^{\alpha}}
\newcommand{\Rnugammad}{R_{\nu,\bsgamma,d}^{\alpha}}
\newcommand{\Aalphagammad}{\calA_{\alpha,\bsgamma,d}}
\newcommand{\Halphagammad}{H_{\alpha,\bsgamma,d}}
\newcommand{\Knud}{K_{\nu,d}^{\alpha}}
\newcommand{\Rnud}{R_{\nu,d}^{\alpha}}
\newcommand{\phih}{\widecheck{\phi}}

\newcommand{\phiKnugamma}{\phi_{\Knugammad}}
\newcommand{\phiRnugamma}{\phi_{\Rnugammad}}
\newcommand{\phiQmgamma}{\phi_{\Qmgammad}}
\newcommand{\half}{\frac{1}{2}}
\DeclareMathOperator{\err}{err}
\DeclareMathOperator{\app}{app}

\newcommand{\dx}{{\rd x}}

\newcommand{\dbx}{{\rd \bsx}}

\newcommand{\gh}{\widehat{g}}

\newcommand{\mubar}{\overline{\mu}}

\DeclareMathOperator{\supp}{supp}

\usepackage{exscale}
\usepackage{epsfig}
\usepackage{verbatim}
\usepackage{a4wide}
\usepackage{graphicx}
\usepackage{multicol}
\usepackage{enumitem}
\usepackage{hyperref}
\usepackage[dvipsnames]{xcolor}
\usepackage{framed}
\usepackage[cal=pxtx]{mathalfa}
\usepackage{bbm}
\usepackage{orcidlink}
\usepackage{soul}
\usepackage{cleveref}

\usepackage[
   doi=false,isbn=false,url=false,eprint=false, 
   maxnames=10,
   giveninits=true,
   style=numeric,
   sorting=nyt]{biblatex}
\setlist[enumerate]{topsep=0pt,itemsep=1ex,partopsep=1ex,parsep=1ex}

\theoremstyle{plain}
\newtheorem{lemma}{Lemma}[section]
\newtheorem{theorem}[lemma]{Theorem}%[section]
\newtheorem{corollary}[lemma]{Corollary}%[section]
\theoremstyle{definition}

\newtheorem{remark}[lemma]{Remark}%[section]
\crefname{lemma}{lemma}{lemmas}
\Crefname{lemma}{Lemma}{Lemmas}
\crefname{theorem}{theorem}{theorems}
\Crefname{theorem}{Theorem}{Theorems}
\crefname{corollary}{corollary}{corollaries}
\Crefname{corollary}{Corollary}{Corollaries}
\crefname{proposition}{proposition}{propositions}
\Crefname{proposition}{Proposition}{Propositions}
\crefname{remark}{remark}{remarks}
\Crefname{remark}{Remark}{Remarks}
\crefname{example}{example}{examples}
\Crefname{example}{Example}{Examples}
\crefname{assumption}{assumption}{assumptions}
\Crefname{assumption}{Assumption}{Assumptions}

\definecolor{darkgreen}{rgb}{0,0.52,0}

\definecolor{darkred}{rgb}{0.72,0,0}

\title[Data Compression using Rank-1 Lattices]{Data Compression using Rank-1 Lattices for Parameter  Estimation in Machine Learning}
\author[M. Gnewuch, K. Harsha, and M. Wnuk]{M. Gnewuch \orcidlink{0000-0002-5516-6476}, K. Harsha \orcidlink{0000-0002-3865-5286}, and M. Wnuk}
\address{Institute for Mathematics, Osnabr\"uck University, Albrechtra{\ss}e 28a, 49074 Osnabr\"uck, Germany}
\email{michael.gnewuch@uni-osnabrueck.de, kumar.harsha@uni-osnabrueck.de, marcin.wnuk@uni-osnabrueck.de}
\date{\today}

\begin{document}

\begin{abstract}
    The mean squared error and regularized versions of it are standard loss functions in supervised machine learning. However, calculating these losses for large data sets can be computationally demanding. Modifying an approach of J. Dick and M. Feischl [Journal of Complexity 67 (2021)], we present algorithms to reduce extensive data sets to a smaller size using rank-1 lattices. 
    Rank-1 lattices are quasi-Monte Carlo (QMC) point sets that are, if carefully chosen, well-distributed in a multidimensional unit cube. 
    The compression strategy in the preprocessing step assigns every lattice point a pair of weights depending on the original data and responses, representing its relative importance. As a result, the compressed data makes iterative loss calculations in optimization steps much faster. 
    We analyze the errors of  our QMC data compression algorithms and the cost of the preprocessing step for functions whose Fourier coefficients decay sufficiently fast so that they lie in certain Wiener algebras or Korobov spaces. In particular, we prove that our approach can lead to arbitrary high convergence rates as long as the functions are sufficiently smooth.
\end{abstract}

\maketitle


\section{Introduction}

Let $\calX = \{ \bsx_1,\dots,\bsx_N\} \subset [0,1]^d$ be a set of data points and $\calY = \{ y_1,\dots,y_N \} \subset \R$ be the corresponding responses. In the context of supervised learning, the goal is to learn a function with a certain structure that maps new data to expected output values. More specifically, we are interested in algorithms that learn parameterized functions $f_\bstheta : [0,1]^d \to \R$ by estimating $\bstheta$ using (regularized) least squares approaches. For such approaches, the quantity of interest is the (regularized) mean squared error 
\begin{equation}\label{term:reg_error}
e(f_\bstheta,\calX,\calY) + \lambda r(\bstheta),
\end{equation}
where
\begin{equation}\label{eq:error_og}
    e(f_\bstheta,\calX,\calY) = \frac{1}{N} \sum_{n=1}^N (f_\bstheta(\bsx_n) - y_n)^2 = \frac{1}{N}\sum_{n=1}^N f_\bstheta^2(\bsx_n) - \frac{2}{N}\sum_{n=1}^N y_n f_\bstheta(\bsx_n) + \frac{1}{N}\sum_{n=1}^N y_n^2,
\end{equation}
$r(\bstheta)$ is a regularization term to penalize model complexity, and $\lambda > 0$ is the regularization strength.

Let us mention some customary choices of regularized errors. 
As usual, we denote the $\ell^p$-norm in $\R^d$ by $\|\cdot \|_p$, i.e., for $\bsx = [x_1,\ldots,x_d] \in \R^d$ we have $\norm{\bsx}_p := \left( \sum_{j \in [d]} \abs{x_j}^p \right)^{1/p}$ if $1\le p <\infty$ and $\norm{\bsx}_{\infty} := \max_{j \in [d]} \abs{x_j}$, where $[d] := \{1,2, \ldots, d\}$. Furthermore, we denote by $\norm{\bsx}_0$ the number of non-zero entries of $\bsx$.
In the case without regularization, i.e., where $r=0$, the loss function in \eqref{term:reg_error} is simply the classical mean squared error as given in \eqref{eq:error_og}, which is the standard loss function in simple linear regression, see e.g. \cite{Theil1997}. 
In the case with regularization the following examples of loss functions are commonly used: 
\begin{enumerate}[label=(\alph*)]
    \item Loss function for \textit{linear regression with the best subset selection} $ e(f_{\bstheta}, \calX, \calY) + \lambda \norm{\bstheta}_0 $, cf. \cite{Natarajan1995,Miller2002}.  
    \item Loss function with \textit{lasso penalty} $ e(f_{\bstheta}, \calX, \calY) + \lambda \norm{\bstheta}_1 $, see \cite{Tibshirani1996}.
    \item Loss function for \textit{ridge regression with Tikhonov regularization} $ e(f_{\bstheta}, \calX, \calY) + \lambda \norm{T \bstheta}_2^2 $, where $T$ is a suitable Tikhonov matrix, cf. \cite{HK_ridge2_1970,HK_ridge_1970}.
    \item Loss function for \textit{elastic net regression} $ e(f_{\bstheta}, \calX, \calY) + \lambda \left( \alpha \norm{\bstheta}_1 + (1-\alpha) \norm{\bstheta}_2^2 \right) $. This approach tries to find a reasonable balance between the advantages of the ridge regression penalty ($\alpha = 0$) and the lasso penalty ($\alpha = 1$), see e.g. \cite{FHT2010}.
\end{enumerate}
Let us emphasize that the data compression approach discussed in this paper works for all regularization functions $r$ in \eqref{term:reg_error}, as long as they do not depend on the data points $\calX$ or the responses $\calY$. 

The goal is to find parameters $\bstheta$ such that the loss function \eqref{term:reg_error} is sufficiently small or even minimal. Now the optimization procedure that is used to calculate suitable parameters $\bstheta$ may rely on a large number of evaluations of the loss function; this applies, e.g., for gradient descent algorithms, as well as for general-purpose optimization using one-shot optimization \cite{BGKTV17}, simulated annealing \cite{Belisle1992}, the Nelder-Mead Simplex algorithm \cite{NM1965}, or the Broyden-Fletcher-Goldfarb-Shanno algorithm \cite{Broyden1970,Fletcher1970,Goldfarb1970,Shanno1970}.
In this case, the first two sums on the right-hand side of \eqref{eq:error_og} lead to cost that is proportional to the number of evaluations of the loss function times the number $N$ of data points and corresponding responses times the cost of the evaluation of a single function $f_\bstheta$. 
Note that the third sum on the right-hand side of \eqref{eq:error_og} does not depend on $\bstheta$ and has to be evaluated at most once, namely if one is interested in the actual values of the errors for different functions $f_\bstheta$.  
If one only wants to find the parameters $\bstheta$ that minimize the error function, then the third sum has not to be evaluated at all, 
since it is independent of $\bstheta$ and can be viewed as an additive constant in the objective function.
Note furthermore, that the regularization term in \eqref{term:reg_error} does not depend on the data. 
Since in present-day applications the amount of data is typically very large, Dick and Feischl proposed in \cite{DF2020} the following approach to alleviate the cost burden
that is caused by the first two sums on the right-hand side of \eqref{eq:error_og}. 
Let $\nu$ be a positive real number that will be specified later. %Following the approach of Dick and Feischl \cite{DF2020}, 
For $L \ll N$, a suitably chosen
point set $\calZ = \{ \bsz_0,\dots,\bsz_{L-1}\}$, and suitably chosen weights $\WXZl,\WXYZl$, we may approximate the mean squared error $e(f_\bstheta,\calX,\calY)$
by the following quantity: 
\begin{equation}\label{eq:error_app}
    \app_L(f_\bstheta,\calX,\calY) := \frac{1}{L}\sum_{\ell=0}^{L-1} f_\bstheta^2(\bsz_\ell) \WXZl - \frac{2}{L}\sum_{\ell=0}^{L-1} f_\bstheta(\bsz_\ell) \WXYZl + \frac{1}{N}\sum_{n=1}^N y_n^2.
\end{equation}
Clearly, evaluating $\app_L(f_\bstheta,\calX,\calY)$ many times for different parameter sets $\bstheta$ would be much cheaper than the corresponding  evaluations of $e(f_\bstheta,\calX,\calY)$.
The main question is now, how to choose the point set $\calZ$ and how to calculate for given data $\calX$ and $\calY$ the corresponding weights $\WXZl$ and  $\WXYZl$?
 Our priority objective is that $\app_L(f_\bstheta,\calX,\calY)$ should be a tight approximation of $e(f_\bstheta,\calX,\calY)$. Another objective is that the cost of computing the weights should be low. Since we can reuse the (pre-)computed weights in the course of the optimization procedure to calculate the loss function for as many sets of parameters $\bstheta$ as we like, this cost can be viewed as precomputation or ``start-up'' cost. Hence, it would be perfectly in order if this cost is actually not very small, but of some reasonable amount instead.

Dick and Feischl recommended choosing $\calZ$ as a digital net, which is a quasi-Monte Carlo point set that is uniformly distributed in the unit cube $[0,1]^d$, cf. \cite{Nie87, DP10}. They showed how to choose the corresponding  weights $\WXZl$ and  $\WXYZl$, estimated the cost of their precomputation, and studied the resulting approximation error 
\begin{equation}\label{total_error}
{\rm err}:= |e(f_\bstheta,\calX,\calY) - \app_L(f_\bstheta,\calX,\calY)|.
\end{equation}
We propose to choose $\calZ$ as a rank-$1$ lattice, constructed with the help of a fast component-by-component algorithm, see \cite{SR2002, Kuo03, NC2006}.
Similarly as in \cite{DF2020}, we derive formulas for the corresponding weights, see Section~\ref{se:derivation}. Furthermore, we estimate the cost of their precomputation and  analyze the approximation error of this approach, see  Section~\ref{se:err_cost}.

We compare both approaches in Section~\ref{se:comparison}. Roughly speaking, the results of the comparison are as follows: On the one hand, we are able to prove for our specific approach 
convergence rates of approximation errors that are better than the ones obtained in \cite{DF2020}. In particular, the (polynomial) convergence rates derived in \cite{DF2020} are always strictly less than $1$, 
regardless of the smoothness of the function $f_\bstheta$, while our approach leads to convergence rates that essentially increase linearly with the smoothness of $f_\bstheta$ and can therefore be arbitrarily large. On the other hand, the cost of precomputing the weights is roughly comparable, but, nevertheless, may be a bit smaller in the case of digital nets. 
As explained above, the computation of the weights can be viewed as precomputation, which is done before the optimization procedure itself starts.
That is why we believe that the (reasonably small) differences in the cost of these computations are not an issue.
Nevertheless, as explained at the end of Section~\ref{se:comparison}, we can easily modify our approach and reduce the cost of precomputing the weights
if we are ready to accept a reduced convergence rate of the approximation error. With this trade-off, we are able to achieve essentially the same cost of precomputation as Dick and Feischl do. For details, see Section~\ref{se:comparison}.

Let us emphasize that the analysis in \cite{DF2020} as well as the analysis in the main part of our paper do not assume that the models, i.e., the functions $f_\bstheta$, have any structure, apart from a certain degree of smoothness, or that the data points $\calX$ have some specific distribution. If the models or data points have additional structure, one can try to exploit this. 
Consider, e.g., the case that the model is 
a truncated Fourier series
$$f_{\bstheta}(\bsx) := \sum_{\bsk \in K} \theta_{\bsk} \omega_{\bsk}(\bsx),$$
    where $K$ is a finite subset of $\Z^d$, $\omega_\bsk$ are the trigonometric basis functions defined as
    \begin{equation*}
        \omega_\bsk(\bsx) := \exp \left(2 \pi i \cdot \bsk \cdot \bsx\right),
    \end{equation*}
    $\bsk \cdot \bsx:= \sum_{j=1}^d k_j x_j$ the standard scalar product in $\R^d$, and $\bstheta := \{ \theta_{\bsk} \}_{\bsk \in K}$ are the Fourier coefficients (parameters) we want to learn. The Fourier coefficients for a function $f \in L_1([0,1]^d)$ are defined as
    \begin{equation*}
        \fh_\bsk := \int_{[0,1]^d} f(\bsx) \exp(-2\pi i \bsk \cdot \bsx) \dbx, \quad \bsk \in \Z^d.
    \end{equation*}
In this case our lattice rule approach may be particularly beneficial, since evaluating the  functions $f_\bstheta$  at all the points of a rank-$1$ lattice $\calZ$ may be done with the help of the univariate fast Fourier transform, regardless of the form of the finite index set $K$, cf. \cite{Kammerer13}. 
This leads to a reasonable speed-up: While the cost of evaluating $f_\bstheta$ in all the points of a  general point set $\calZ$ is of order $O(d|K|L)$, choosing $\calZ$ as a rank$1$-lattice reduces it to $O(d|K| + L \log(L))$, cf. Section~\ref{se:precomp_gen_K}.
Truncated Fourier series models are of particular interest in the explainable ANOVA approximation approach in machine learning as introduced in \cite{PS2021a}  and further developed, e.g., in \cite{PS2021, PS2022}.

Of course, there are known approaches that are different from, but related to, those proposed by Dick and Feischl and in this paper such as subsampling with or without leveraging, coresets or support points; for a discussion of relevant algorithms, their convergence rates and cost, as well as references to the literature, we refer the reader to Section~1.1 of \cite{DF2020} and to \cite{FSS13}, \cite{WZM18}.

In the article we are considering many different cases, so let us give here further information about its
structure and an overview 
of the different settings we analyze. 
In Section 2 we derive the compression weights. To this end, we  expand the functions $f_{\theta}$ and $f_{\theta}^2$,
that appear in the first two sums in  \eqref{eq:error_og}, into a Fourier series and then truncate it. More precisely, 
for $g\in \{f_{\theta}, f_{\theta}^2\}$ we choose a finite set of frequencies $K_\nu \subset \Z^d$ and use the truncated Fourier series
$g_{K_\nu} := \sum_{\bsk \in K_\nu} \widehat{g}_{\bsk} \omega_{\bsk}$ to approximate $g$. The compression weights and their computing time depend 
crucially on the  choice of the set of frequencies $K_\nu$, and the size of the total error ${\rm err}$ in \eqref{total_error} does as well.
For our error analysis in Section~\ref{se:err_cost} 
we assume that the functions $g\in \{f_{\theta}, f_{\theta}^2\}$ belong to
\begin{itemize}
\item weighted Korobov spaces $\Halphagammad$  
\item or weighted Wiener algebras $\Aalphagammad$.
\end{itemize}
The error analysis relies on the fact that the total error ${\rm err}$ in \eqref{total_error} can be decomposed 
into a truncation error ${\rm err}_1$ (induced by the choice of the frequency set $K_\nu$), cf. \eqref{eq:error1},
and an approximation error ${\rm err}_2$, cf. \eqref{eq:error2}.
Due to  a suitable  multiplier estimate the latter error can be estimated with the help of a bound for the integration 
error of lattice rules. Consequently, our error analysis
consists of three steps:
\begin{itemize}
\item bounding the truncation error ${\rm err_1}$, 
\item providing suitable multiplier estimates,
\item and using the right bound for the integration error of lattice rules.
\end{itemize}
Those steps are described - without appealing to concrete choices of the set of frequencies $K_\nu$ - in 
Sections~\ref{SubSec_Trunc_Err}
to \ref{SubSec_Lattice}. In Section~\ref{se:precomp_gen_K} we discuss the cost of computing the compression weights for a general set $K_\nu$. 
Afterwards we move to the error analysis for concrete choices of $K_\nu$, namely
\begin{itemize}
\item the weighted continuous hyperbolic cross $\Knugammad$,
\item the weighted $d$-dimensional rectangle $\Rnugammad$,
\item and the weighted step hyperbolic cross $\Qmgammad$.
\end{itemize}
These three choices are discussed in Sections~\ref{se:contX} to \ref{se:stepX}.
In Sections~\ref{se:precomp_wtH}, 
\ref{Subsec:Rectangle_Precomputation}, 
and \ref{subsec:precomp_stepX}
we describe how the compression weights corresponding to these three different choices of $K_\nu$ can be precomputed efficiently.
As already mentioned before, we compare our results to the ones from \cite{DF2020} in Section~\ref{se:comparison}; 
in particular, we compare the precomputation times for both approaches in Section~\ref{Subsec:Comparison_Cost}.

Let us finish this introduction by emphasizing that this paper focuses solely on the theoretical analysis of the QMC data compression approach based on lattice point sets.
In particular, our main goals are to provide a transparent error analysis that may also be extended to other interesting point sets, and to show that the convergence rates of approximation errors can indeed be arbitrarily high if the input functions are sufficiently smooth. A rigorous study of numerical experiments for interesting use cases is beyond the scope of this paper.

\section{Derivation of compression weights}\label{se:derivation}

In \cite{DF2020} the authors presented two approaches to deduce formulas for the compression weights from \eqref{eq:error_app} that are suitable for digital nets $\calZ$.
The first approach is based on geometrical properties of digital nets and on an inclusion-exclusion formula. It seems to us that this approach cannot be modified in a reasonable way to work also for lattice point sets. The second approach is based on a truncation of the Walsh series of the functions $f_\bstheta$ and $(f_\bstheta)^2$, respectively.
As we will show now, this approach can be modified appropriately for lattice point sets $\calZ$ by considering Fourier series instead of Walsh series. As in \autocite{DF2020}, we consider a function $g$ and coefficients $c_n$, where, in order to derive the weights $\WXZl$ we set
\begin{equation}\label{eq:gcn1}
    g = f_\bstheta^2, \quad c_n = 1,
\end{equation}
and for the weights $\WXYZl$, we set
\begin{equation}\label{eq:gcn2}
    g = f_\bstheta, \quad c_n = y_n.
\end{equation}
Let us assume that the Fourier series of $g$ converges absolutely and uniformly on $[0,1]^d$.
Let $K \subset \Z^d$ be a finite subset that will be suitably chosen later. Then we define the function
\begin{equation*}
        g_K := \sum_{\bsk \in K} \ghk \omega_{\bsk}, 
\end{equation*}
where $\ghk = \innProd{g}{\omega_{\bsk}}{L^2([0,1]^d)}$ is the $\bsk$-th Fourier coefficient of $g$. 
Consequently,
\begin{equation*}
        g - g_{K} = \sum_{\bsk \in \Z^d \setminus K} \ghk \omega_{\bsk}.
\end{equation*}
We want to choose $K$ big enough so that the tail $g-g_{K}$ is small. In other words, we aim to have $g_K \approx g$.

Assume that $K = K_{\nu,d}$ depends on a parameter $\nu > 0$ such that $\lim_{\nu \to \infty} \norm{g - g_{K_{\nu,d}}}_{\infty} = 0.$ The parameter $\nu$ is allowed to depend on $d$. Then, from the orthogonality of the functions $\omega_{\bsk}$ in $L^2([0,1]^d)$ we have
\begin{equation}\label{eq:approx1}
    \begin{aligned}
        \frac{1}{N} \sum_{n=1}^N c_n g(\bsx_n) %&= \frac{1}{N} \sum_{n=1}^N c_n \sum_{\bsk \in \Z^d} \left( \int_{[0,1]^d} g(\bsx) \overline{\omega_{\bsk}(\bsx)} \dbx \right) \omega_{\bsk}(\bsx_n) \\
        &= \frac{1}{N} \sum_{n=1}^N c_n g_K(\bsx_n) + \frac{1}{N} \sum_{n=1}^N c_n (g-g_{K})(\bsx_n)\\
        &\approx \frac{1}{N} \sum_{n=1}^N c_n g_{K}(\bsx_n) = \frac{1}{N} \sum_{n=1}^N c_n \sum_{\bsk \in K} \ghk \omega_{\bsk}(\bsx_n)\\
        &= \int_{[0,1]^d} g(\bsx) \left( \sum_{\bsk \in K} \overline{\omega_{\bsk}(\bsx)} \frac{1}{N} \sum_{n=1}^N c_n \omega_{\bsk}(\bsx_n) \right) \dbx \\
        &= \int_{[0,1]^d} g(\bsx) \phi_{K}(\bsx) \dbx,
    \end{aligned}
\end{equation}
for the function
$$ \phi_K(\bsx) := \sum_{\bsk \in K} \phih_{\bsk} \overline{\omega_{\bsk}(\bsx)} $$
where the coefficients are
\begin{equation*}
    \phih_{\bsk} := \frac{1}{N} \sum_{n=1}^N c_n \omega_{\bsk}(\bsx_n).
\end{equation*}

Let us write $\calC := ( c_n)_{n \in [N]}$. We define the truncation error, obtained from \eqref{eq:approx1}, as
\begin{equation}\label{eq:error1}
    \begin{aligned}
        \err_1(g,\calC) &:= \abs{\frac{1}{N} \sum_{n=1}^N c_n g(\bsx_n) - \int_{[0,1]^d} g(\bsx) \phi_{K}(\bsx) \dbx} = \abs{\frac{1}{N} \sum_{n=1}^N c_n (g-g_{K})(\bsx_n)} \\
        &\leq \norm{g-g_{K}}_{\infty} \frac{1}{N} \sum_{n=1}^N \abs{c_n}.
    \end{aligned}
\end{equation}

Once $K$ is chosen, we can approximate the last integral in \eqref{eq:approx1} by a suitable quadrature given by points $\calZ = \{ \bsz_0,\dots,\bsz_{L-1}\}$ as follows
\begin{equation}\label{eq:approx2}
    \begin{aligned}
        \int_{[0,1]^d} (g \phi_{K})(\bsx) \dbx \approx \frac{1}{L} \sum_{\ell = 0}^{L-1} (g \phi_{K})(\bsz_{\ell}).
    \end{aligned}
\end{equation}
Consequently, we define the approximation error obtained from \eqref{eq:approx2} as
\begin{equation}\label{eq:error2}
    \begin{aligned}
        \err_2(g,\calC) := \abs{\int_{[0,1]^d} (g \phi_{K})(\bsx) \dbx - \frac{1}{L} \sum_{\ell = 0}^{L-1} (g \phi_{K})(\bsz_{\ell})}.
    \end{aligned}
\end{equation}
Both steps \eqref{eq:approx1} and \eqref{eq:approx2} lead to the approximation
\begin{equation*}
    \frac{1}{N} \sum_{n=1}^N c_n g(\bsx_n) \approx \int_{[0,1]^d} g(\bsx) \phi_{K}(\bsx) \dbx \approx \frac{1}{L} \sum_{\ell = 0}^{L-1} (g \phi_{K})(\bsz_{\ell}).
\end{equation*}
If we follow the approach described above, the weights $\WXZl$ and $\WXYZl$ in \eqref{eq:error_app} are of the form
\begin{equation}\label{eq:phiK}
    \phi_K(\bsz_\ell) = \frac{1}{N} \sum_{\bsk \in K} \sum_{n=1}^N c_n \omega_{\bsk}(\bsx_n - \bsz_\ell).
\end{equation}
More precisely, for $g = f_\bstheta^2$ and $c_n = 1$, we get the first set of input-dependent weights
$$ W_{\calX,\calZ,\nu, \ell} = \frac{1}{N}\sum_{\bsk \in K} \sum_{n=1}^N \omega_{\bsk}(\bsx_n - \bsz_\ell). $$
For $g = f_\bstheta$ and $c_n = y_n$, we get the second set of input-output dependent weights
$$ W_{\calX,\calY,\calZ,\nu, \ell} = \frac{1}{N} \sum_{\bsk \in K} \sum_{n=1}^N y_n \omega_{\bsk}(\bsx_n - \bsz_\ell). $$

\begin{remark}\label{re:err}
    The total error we are interested in is
    \begin{equation*}
        \begin{aligned}
            \abs{e(f_\bstheta,\calX,\calY) - \app_L(f_\bstheta,\calX,\calY)} \leq & \abs{\frac{1}{N}\sum_{n=1}^N f_\bstheta^2(\bsx_n) - \frac{1}{L}\sum_{\ell=0}^{L-1} f_\bstheta^2(\bsz_\ell) \WXZl}\\ & + 2\abs{\frac{1}{N}\sum_{n=1}^N y_n f_\bstheta(\bsx_n) - \frac{1}{L}\sum_{\ell=0}^{L-1} f_\bstheta(\bsz_\ell) \WXYZl}\\
            &\leq \err_1(f_\bstheta^2, \bsone_N) + \err_2(f_\bstheta^2, \bsone_N) + 2(\err_1(f_\bstheta, \calY) + \err_2(f_\bstheta, \calY)),
        \end{aligned}
    \end{equation*}
    where $\bsone_N:= (1)_{n\in[N]}$.
    That is, to estimate the error, it suffices to estimate
    \begin{equation}\label{eq:error_gcn}
        \abs{\frac{1}{N} \sum_{n=1}^N c_n g(\bsx_n) - \frac{1}{L} \sum_{\ell = 0}^{L-1} (g \phi_{K})(\bsz_{\ell})} \leq \err_1(g,\calC) + \err_2(g,\calC),
    \end{equation}
    which follows from the triangle inequality for $g$ and $c_n$ as in \RefEqTwo{eq:gcn1}{eq:gcn2}, $\err_1$ as in \eqref{eq:error1}, $\err_2$ as in \eqref{eq:error2}, and $\phi_K$ as in \RefEq{eq:phiK}.
\end{remark}

\begin{remark}
In this paper we focus on approximating the error for  periodic functions $f$. The uniformly and absolutely convergent expansion of $f$ allows us to perform a suitable truncation to the set of the presumably most important frequencies and the lattice rules which we use to approximate the integral work well with periodic functions. Similar strategies (possibly with different series expansion and integration rules) could be used to tackle also  non-periodic functions $f$ (as already done in \cite{DF2020}), but this is beyond the scope of the present article.
\end{remark}

\section{Error and cost analysis}\label{se:err_cost}

For our error and cost analysis, we consider weighted function spaces with product weights induced by a decreasing sequence of weights 
$$1 \geq \gamma_1 \geq \gamma_2 \geq ... \geq \gamma_d > 0$$ 
assigned to all the variables. Then the weight (or importance) of a set of interacting variables $\bsu$ is given by $\gamma_\bsu = \prod_{j\in \bsu} \gamma_j$ for non-empty $\bsu$, and we set $\gamma_{\emptyset} = 1$. Product weights for weighted function spaces were first introduced in \cite{SW1998}.

Let $\alpha > 1/2$ and define 
$$r_{\alpha}(\bsgamma, \bsh) := \prod_{j \in [d]} r_\alpha(\gamma_j,h_j), \quad \text{ where }\quad r_\alpha(\gamma_j,h_j) = \max \left( \frac{\abs{h_j}^{2\alpha}}{\gamma_j}, 1 \right).$$

The \textbf{weighted Wiener algebra} $\Aalphagammad$, see \cite{BD73} and \cite{KLT2023}, is defined as
\begin{equation}\label{eq:wtW}
    \Aalphagammad := \left\{ f \in L_1([0,1]^d) : \norm{f}_{\Aalphagammad} := \sum_{\bsk \in \Z^d} \sqrt{r_{\alpha}(\bsgamma, \bsk)} \abs{\fh_{\bsk}} < \infty \right\}.
\end{equation}

Even though we only enforce $L_1$-integrability to define the Wiener algebra in \eqref{eq:wtW}, its elements are nevertheless $L_2$-integrable. This can be shown using the Plancherel theorem as follows
\begin{equation*}
    \norm{f}_{L_2([0,1]^d)}^2 = \sum_{\bsk \in \Z^d} \abs{\fh_\bsk}^2 \leq \max_{\bsell \in \Z^d} \frac{1}{r_\alpha(\bsgamma,\bsell)} \sum_{\bsk \in \Z^d} r_\alpha(\bsgamma,\bsk) \abs{\fh_\bsk}^2 \leq \left(\sum_{\bsk \in \Z^d} \sqrt{r_\alpha(\bsgamma,\bsk)} \abs{\fh_\bsk} \right)^2 = \norm{f}_{\Aalphagammad}^2.
\end{equation*}

The \textbf{weighted Korobov space} $\Halphagammad$ is defined as
\begin{equation}\label{eq:wtH}
    \Halphagammad := \left\{ f \in L_2([0,1]^d) : \norm{f}_{\Halphagammad} := \sqrt{\sum_{\bsk \in \Z^d} r_\alpha(\bsgamma,\bsk) \abs{\fh_{\bsk}}^2} < \infty \right\}.
\end{equation}
Clearly, as vector spaces, $\Aalphagammad$ and $\Halphagammad$ do not depend on the weights $\gamma_j$.
Note that $\Aalphagammad \subset \Halphagammad$. The smoothness parameter $\alpha$ controls the rate of decay of the Fourier coefficients of the functions and guarantees that the functions exhibit a certain degree of smoothness. 
For $\alpha > 1/2$, the Fourier series of any function $f \in \Halphagammad$ converges absolutely and uniformly on $[0,1]^d$ since
\begin{equation*}
    \sum_{\bsk \in \Z^d} \abs{\fh_{\bsk}} \leq \norm{f}_{\Halphagammad} \cdot \left( \sum_{\bsk \in \Z^d} \frac{1}{r_{\alpha}(\bsgamma,\bsk)} \right)^{1/2}  < \infty,
\end{equation*}
where the first inequality follows from the Cauchy-Schwarz inequality. Consequently, all elements of $\Halphagammad$ are continuous periodic functions.
If $\alpha$ is an integer, any element in the space $H_{\alpha,\gamma,1}$ of univariate functions has $\alpha$ weak derivatives, where the $k$-th derivative is absolutely continuous for $k = 1,2,\dots,\alpha-1$, and the $\alpha$-th derivative is square integrable, 
and the norm of $H_{\alpha,\gamma,1}$
fulfills 
\begin{equation}\label{korobov_norm_derivatives}
   \norm{f}_{H_{\alpha,\gamma,1}}^2 
   = \left| \int_0^1 f(x)\, dx \right|^2 + \frac{1}{\gamma} \frac{1}{(2\pi)^{2\alpha}} \int_0^1 |f^{(\alpha)} (x) |^2 \,dx,
\end{equation}
see, e.g., \cite[Appendix~A.1]{NW08}. Moreover, elements of the spaces $\Halphagammad$ of multivariate functions have weak mixed derivatives up to order $\alpha$ in each variable. Conversely: let $\alpha \geq 1$ be such that $2\alpha \in \mathbb{N}$. If a $1-$periodic function $f$ has all the continuous mixed derivatives $D^{\beta}f$ for $\beta \in \{0, \ldots, 2\alpha  \}^d$ then it is in  $\Halphagammad$. For details see Prop 4.16. in \cite{LP14} (note that the authors parametrize the Korobov spaces differently.)

\begin{remark}
We use the product structure of the weights in a threefold way. Firstly, it gives us the inequality (\ref{eq:r_ineq_wt}) needed in the analysis of multiplier estimates. Secondly, we rely on the results on numerical integration for  weighted Korobov spaces with  product weights, see (\ref{eq:lattice_wt_err}). Finally, the product structure plays a role when analysing the concrete choice of $K$ (weighted rectangle and both types of hyperbolic crosses). We believe that with some effort our results could be generalized to, e.g., product order-dependent weights. However, obtaining similar results for general weights seems to require more work.   
\end{remark}

\subsection{Truncation errors}\label{SubSec_Trunc_Err}

In this section, we provide general estimates for the truncation error $\err_1(g,\calC)$, defined as in \eqref{eq:error1}. More precisely, we rely on the estimate $$\err_1(g,\calC) \leq \norm{g - g_K}_{\infty} \cdot \frac{1}{N} \sum_{n=1}^N \abs{c_n}$$ from \eqref{eq:error1}, and therefore we only need to take care of the term $\norm{g - g_K}_{\infty}$.

\begin{lemma}\label{le:trunc_wtW}
    Let $g \in \Aalphagammad$ and $K \subset \Z^d$ be finite. Then we have that the truncation error satisfies
    $$ \norm{g - g_{K}}_{\infty} \leq \norm{g}_{\Aalphagammad} \left( \sup_{\bsk \notin K} \frac{1}{r_{\alpha}(\bsgamma, \bsk)} \right)^{1/2}. $$
\end{lemma}
\begin{proof}
    We use the triangle inequality
    \begin{equation*}
        \begin{aligned}
            \norm{g-g_{K}}_{\infty} &\leq \sum_{\bsk \notin K} \abs{\gh_{\bsk}} = \sum_{\bsk \notin K} \abs{\gh_{\bsk}} \sqrt{\frac{r_\alpha(\bsgamma,\bsk)}{r_\alpha(\bsgamma,\bsk)}} \\
            &\leq \sup_{\bsk \notin K} r_{\alpha}(\bsgamma, \bsk)^{-1/2} \sum_{\bsk \notin K} \abs{\ghk} \sqrt{r_{\alpha}(\bsgamma,\bsk)} \\
            &\leq  \norm{g}_{\Aalphagammad} \sup_{\bsk \notin K} r_{\alpha}(\bsgamma, \bsk)^{-1/2},
        \end{aligned}
    \end{equation*}
    and the result follows.
\end{proof}

\begin{lemma}\label{le:trunc_wtH}
    Let $g \in \Halphagammad$ and $K \subset \Z^d$ be finite. Then we have that the truncation error satisfies
    $$ \norm{g - g_{K}}_{\infty} \leq \norm{g}_{\Halphagammad} \left( \sum_{\bsk \notin K} \frac{1}{r_{\alpha}(\bsgamma, \bsk)} \right)^{1/2}. $$
\end{lemma}
\begin{proof}
    We use the triangle inequality and the Cauchy-Schwarz inequality
    \begin{equation*}
        \begin{aligned}
            \norm{g-g_{K}}_{\infty} &\leq \sum_{\bsk \notin K} \abs{\gh_{\bsk}} = \sum_{\bsk \notin K} \abs{\gh_{\bsk}} \sqrt{\frac{r_\alpha(\bsgamma,\bsk)}{r_\alpha(\bsgamma,\bsk)}} \\
            &\leq \left( \sum_{\bsk \notin K} \abs{\ghk}^2 r_\alpha(\bsgamma,\bsk) \right)^{1/2} \left( \sum_{\bsk \notin K} 1/r_{\alpha}(\bsgamma,\bsk) \right)^{1/2} \\
            &\leq \norm{g}_{\Halphagammad} \left( \sum_{\bsk \notin K} 1/r_{\alpha}(\bsgamma,\bsk) \right)^{1/2},
        \end{aligned}
    \end{equation*}
    hence the result follows.
\end{proof}

\begin{remark}
    Similar results on the truncation error for specific choices of $K$ have been shown in other papers, e.g. in \cite[Lemma~10]{DF2020},\cite{KPV2015}.
\end{remark}

\subsection{Multiplier estimates}

We notice in \eqref{eq:error2} that $\err_2(g, \calC)$ is actually an integration error implied by an equal-weight cubature rule with respect to the integrand $g\phi_K$, instead of $g$. We exploit this observation in our error analysis and bound the norms of $g\phi_K$ that are of interest to us from above.

Let $\delta > 0$ and put
\begin{equation}\label{eq:const}
    \mubar_N := N^{-1} \sum_{n=1}^N \abs{c_n}, \quad c_{\alpha,\bsgamma,d} := \sqrt{\prod_{j \in [d]} \max(1,2^{2\alpha}\gamma_j)}, \quad \text{ and } \zeta_{\delta,d} := \left[ 1 + 2\zeta(1+2\delta) \right]^{d/2},
\end{equation}
where $\zeta$ is the Riemann zeta function.
Let us state the following inequality shown in \cite{NSW04},
\begin{equation}\label{eq:r_ineq_wt}
    \frac{r_{\alpha}(\bsgamma,\bsell)}{r_{\alpha}(\bsgamma,\bsell + \bsk)} \leq r_{\alpha}(\bsgamma,\bsk) c_{\alpha,\bsgamma,d}^2 \quad \forall \bsk,\bsell \in \Z^d,
\end{equation}
which will be helpful to prove the next two lemmas.

\begin{lemma}\label{le:gphi_wtW}
    Let $\alpha > 1$. If $g \in \Aalphagammad$, then $g \phi_K \in \Halphaminusgamma$ for every $\delta \in (0,\alpha-1)$. Moreover,
    \begin{equation*}
        \norm{g \phi_{K}}_{\Halphaminusgamma} \leq \norm{g}_{\Aalphagammad} \left( \max_{\bsk \in K} r_{\alpha}(\bsgamma,\bsk) \right)^{1/2} c_{\alpha,\bsgamma,d}\ \zeta_{\delta,d}\  \mubar_N,
    \end{equation*}
    where $c_{\alpha,\bsgamma,d}$, $\zeta_{\delta,d}$ and $\mubar_N$ are as defined in \eqref{eq:const}.
\end{lemma}
\begin{proof}
    According to the definition of $\phi_K$, the function $g \phi_K$ can be expressed as
    \begin{equation*}
        \left(g \phi_K\right)(\bsz) = \sum_{\bsk \in K} g(\bsz) \overline{\omega_\bsk(\bsz)} \frac{1}{N} \sum_{n=1}^N c_n \omega_{\bsk}(\bsx_n).
    \end{equation*}
    Consequently, the Fourier coefficients of $g \phi_K$ for $\bsell \in \Z^d$ can be calculated as follows
    \begin{equation}\label{eq:beta_ell}
        \begin{aligned}
            \beta_\bsell := \widehat{\left( g \phi_K \right)}_\bsell &= \int_{[0,1]^d} \left( g \phi_K \right)(\bsz) \overline{\omega_\bsell(\bsz)} \rd \bsz = \sum_{\bsk \in K} \int_{[0,1]^d} g(\bsz) \overline{\omega_\bsk(\bsz)} \overline{\omega_\bsell(\bsz)} \rd \bsz \frac{1}{N} \sum_{n=1}^N c_n \omega_{\bsk}(\bsx_n) \\
            &= \sum_{\bsk \in K} \gh_{\bsell + \bsk} \frac{1}{N} \sum_{n=1}^N c_n \omega_{\bsk}(\bsx_n).
        \end{aligned}
    \end{equation}
    We can now estimate these coefficients using the triangle inequality and \eqref{eq:r_ineq_wt} as follows
    \begin{equation}\label{eq:abs_beta_ell}
        \begin{aligned}
            \abs{\beta_{\bsell}} \leq \mubar_N \sum_{\bsk \in K} \abs{ \gh_{\bsell + \bsk} } &\leq \mubar_N \left( \sum_{\bsk \in K} \abs{ \gh_{\bsell + \bsk} } \sqrt{r_{\alpha}(\bsgamma,\bsell + \bsk)} \right) \cdot \max_{\bsk \in K} \frac{1}{\sqrt{r_{\alpha}(\bsgamma,\bsell + \bsk)}}  \\
            &\leq \norm{g}_{\Aalphagammad} \mubar_N \max_{\bsk \in K} \frac{1}{\sqrt{r_{\alpha}(\bsgamma,\bsell + \bsk)}} \\
            &\leq \norm{g}_{\Aalphagammad} \frac{\mubar_N c_{\alpha,\bsgamma,d}}{\sqrt{r_{\alpha}(\bsgamma,\bsell)}} \max_{\bsk \in K} \sqrt{r_{\alpha}(\bsgamma,\bsk)}.
        \end{aligned}
    \end{equation}
    Using the above estimate, we have
    \begin{equation}\label{eq:wt_gphi}
        \begin{aligned}
            \norm{g \phi_{K}}_{\Halphaminusgamma} &= \left[ \sum_{\bsell \in \Z^d} \abs{\beta_\bsell}^2 r_{\alpha - \half - \delta}(\bsgamma,\bsell) \right]^{1/2}\\
            &\leq \norm{g}_{\Aalphagammad} \mubar_N c_{\alpha,\bsgamma,d} \left[ \sum_{\bsell \in \Z^d} \frac{r_{\alpha-\half-\delta}(\bsgamma,\bsell)}{r_{\alpha}(\bsgamma,\bsell)} \right]^{1/2} \max_{\bsk \in K} \sqrt{r_{\alpha}(\bsgamma,\bsk)} \\
            &= \norm{g}_{\Aalphagammad} \mubar_N\ c_{\alpha,\bsgamma,d}\ \zeta_{\delta,d} \max_{\bsk \in K} \sqrt{r_{\alpha}(\bsgamma,\bsk)},
        \end{aligned}
    \end{equation}
    and the statement follows.
\end{proof}

\begin{lemma}\label{le:gphi_wtH}
    Let $\alpha > 1$. If $g \in \Halphagammad$, then $g \phi_K \in \Halphaminusgamma$ for every $\delta \in (0,\alpha-1)$. Moreover,
    \begin{equation*}
        \norm{g \phi_K}_{\Halphaminusgamma} \leq \norm{g}_{\Halphagammad} \left( \sum_{\bsk \in K} r_{\alpha}(\bsgamma,\bsk) \right)^{1/2} c_{\alpha,\bsgamma,d}\ \zeta_{\delta,d}\  \mubar_N,
    \end{equation*}
    where $c_{\alpha,\bsgamma,d}$, $\zeta_{\delta,d}$ and $\mubar_N$ are as defined in \eqref{eq:const}.
\end{lemma}
\begin{proof}
    We use the Fourier coefficients $\beta_{\bsell}$ for $\bsell \in \Z^d$ of $g\phi_K$ as in \eqref{eq:beta_ell} with the Cauchy-Schwarz inequality and \eqref{eq:r_ineq_wt} to get
    \begin{equation*}
        \begin{aligned}
            \abs{\beta_{\bsell}} \leq \mubar_N \sum_{\bsk \in K} \abs{\gh_{\bsell + \bsk}} &\leq \mubar_N \left( \sum_{\bsk \in K} \abs{\gh_{\bsell + \bsk}}^2 r_{\alpha}(\bsgamma,\bsell + \bsk) \right)^{\half} \left( \sum_{\bsk \in K} \frac{1}{r_{\alpha}(\bsgamma,\bsell + \bsk)} \right)^{\half} \\
            &\leq \norm{g}_{\Halphagammad} \mubar_N \left( \sum_{\bsk \in K} \frac{1}{r_{\alpha}(\bsgamma,\bsell + \bsk)} \right)^{\half} \\ &\leq \norm{g}_{\Halphagammad} \frac{\mubar_N c_{\alpha,\bsgamma,d}}{\sqrt{r_{\alpha}(\bsgamma,\bsell)}} \left( \sum_{\bsk \in K} r_{\alpha}(\bsgamma,\bsk) \right)^{\half},
        \end{aligned}
    \end{equation*}
    where $\mubar_N$ is defined as in \eqref{eq:const}.
    We now proceed similarly as in \eqref{eq:wt_gphi} to get
    \begin{equation*}
        \norm{g\phi_K}_{\Halphaminusgamma} \leq \norm{g}_{\Halphagammad} \mubar_N\ c_{\alpha,\bsgamma,d}\ \zeta_{\delta,d} \left( \sum_{\bsk \in K} r_{\alpha}(\bsgamma,\bsk) \right)^{\half},
    \end{equation*}
    and the statement is proved.
\end{proof}

\subsection{Rank-\texorpdfstring{$1$}{1} Lattice point sets}\label{SubSec_Lattice}

As mentioned in the introduction, to compress large data sets consisting of data points $\calX$ and the corresponding responses $\calY$, we propose to use rank-$1$ lattice point sets $\calZ$. 
On the one hand, lattice rules, i.e., equal-weight cubatures whose integration knots form a lattice point set, can be used to approximate the integral of $g \phi_{K}$ well. This step will enable us to estimate the approximation error $\err_2(g,\calC)$ defined in \eqref{eq:error2}.
On the other hand, we may exploit the structure of lattice rules for the precomputation of the compression weights by making use of the univariate fast Fourier transform or the ``Dirichlet kernel trick'', cf. Sections~\ref{se:precomp_gen_K},
\ref{Subsec:Rectangle_Precomputation}, and \ref{subsec:precomp_stepX}.

\subsubsection{Rank-\texorpdfstring{$1$}{1} Lattice Rules for Multivariate Integration}

For $L\in \N$ and a generating vector $\mathbb{\frakg} \in \{ 1,\dots,L-1 \}^d$, the corresponding lattice point set is given by
\begin{equation}\label{def:Z}
    \calZ := \left\{ \bsz_\ell = \frac{\ell}{L}\mathbb{\frakg} \mod 1 : \ell = 0,1,\dots,L-1 \right\},
\end{equation}
where $\bmod\ 1$ (i.e., taking the fractional part) is meant component-wise.
We mainly focus on the case where $L$ is a prime number. 
We define the worst-case integration error in the Korobov space $\Halphagammad$ as
\begin{equation*}
    e(\Halphagammad, \calZ) = \sup_{\substack{ f \in \Halphagammad,\\ \norm{f}_{\Halphagammad} \leq 1 }} \abs{ \int_{[0,1]^d} f(\bsx) \dbx - \frac{1}{L} \sum_{\ell = 0}^{L-1} f(\bsz_\ell) },
\end{equation*}
which has a closed-form expression, see, e.g., \cite[Theorem~4]{Kuo03} or \cite{DKP2022}, given by
\begin{equation*}
    e(\Halphagammad, \calZ) = -1 + \frac{1}{L} \sum_{\ell=0}^{L-1} \prod_{j \in [d]} \left( 1 + \gamma_j \varphi_\alpha\left(\frac{\ell \frakg_j}{L} \right) \right), \quad \text{where } \varphi_\alpha(x) = \sum_{h \in \Z \setminus \{0\}} \frac{\exp(2\pi i h x)}{\abs{h}^{2\alpha}},
\end{equation*}
and $\frakg_j$ are the components of the generating vector $\mathbb{\frakg}$. If $\alpha$ is an integer, we can express $\varphi_\alpha$ as the Bernoulli polynomial $B_{2\alpha}$ of degree $2\alpha$. The component-by-component(CBC) algorithm constructs the generating vector $\mathbb{\frakg} = (\frakg_j)_{j \in [d]}$ by putting $\frakg_1:=1$, and successively determining the other components, 
where for fixed $\frakg_1, \ldots, \frakg_{j-1}$ the component $\frakg_j$ is obtained by greedily minimizing the closed-form expression of the worst-case integration error in $H_{\alpha,\bsgamma,j}$. 
Using the resulting generating vector, the error satisfies
\begin{equation}\label{eq:lattice_wt_err}
    e(H_{\alpha,\bsgamma,d}, \calZ) \leq C_{\bsgamma,d}(\alpha,\tau) L^{-\alpha+\tau} \quad \forall \tau \in \left( 0, \alpha-\half \right],
\end{equation}
where
\begin{equation*}
    C_{\bsgamma,d}(\alpha,\tau) = 2^{\alpha-\tau} \prod_{j \in [d]} \left[ 1 + 2 \gamma_j^{\frac{1}{2(\alpha-\tau)}} \zeta\left( \frac{\alpha}{\alpha-\tau} \right) \right]^{\alpha - \tau},
\end{equation*}
see \cite[Theorem~4]{Kuo03} and \cite[Chapter~3]{DKP2022}.

\begin{comment}
We also address the unweighted function space setting, i.e. $\gamma_1 = \dots = \gamma_d = 1$, since we will use it later in \RefSec{se:comparison}. If we construct a generating vector for a lattice $\calZ$ of size $L$ via the CBC algorithm, then the worst-case integration error of the corresponding lattice rule satisfies
\begin{equation}\label{eq:lattice_err}
    e(\Halphad, \calZ) \leq (1 + 2\zeta(2\alpha))^{d/2} L^{-\alpha} \left[ 1 + 2^{2\alpha(d+1)}(1 + \log L)^{2\alpha d} \right]^{1/2} \lesssim \left( \frac{\log(L)^d}{L} \right)^{\alpha},
\end{equation}
where $\Halphad$ denotes the unweighted Korobov space $H_{\alpha,\bsone_d,d}$, see e.g. \autocite[Theorem~4.21]{LP14}.    
\end{comment}

We provide a short review of algorithms to construct generating vectors for rank-$1$ lattice rules.
The fast CBC algorithm was introduced in \cite{NC2006}, where the authors reformulated the optimization steps in the CBC algorithm as a circulant matrix-vector product. The circulant matrix under consideration satisfies specific symmetry properties that allow the use of the fast Fourier transform (FFT) to compute the matrix-vector product with a cost of order $d L \log(L)$, in contrast to the cost of order $dL^2$ for the standard CBC algorithm.

The reduced CBC algorithm is another improvement over the standard CBC algorithm to construct generating vectors for the integration problem on $\Halphagammad$. The underlying idea behind the reduced algorithm is to shrink the search spaces for the components $\frakg_j$ of the lattice generating vector $\mathbb{\frakg}$ for which the corresponding weights $\gamma_j$ are small. This idea can be combined with the fast CBC algorithm to form the reduced fast CBC algorithm, see \cite{DKLP2015}.

The generating vector produced by the CBC algorithm can be further refined using the successive coordinate search (SCS) algorithm, introduced in \cite{ELN2018}. The SCS algorithm iteratively improves the provided generating vector by optimizing the error in a component-wise fashion. Since this approach is similar to the CBC algorithm, the SCS algorithm can be accelerated using its reduced fast version; see \cite{EK2019}.

In summary, there are efficient algorithms available to construct the generating vector of a lattice rule that guarantees a small integration error on weighted Korobov spaces. For a comprehensive overview of the CBC algorithm and its variants, we refer to \cite[Chapters~3,4]{DKP2022}.

\subsubsection{Precomputation of the compression weights for a general index set $K$}\label{se:precomp_gen_K}

We now discuss the precomputation of the compression weights for a rank-$1$-lattice point set $\calZ$ as in \eqref{def:Z} and an arbitrary finite set $K\subset \Z^d$. 
We present a general approach to %efficiently
compute the weights $\phi_K(\bsz_\ell)$, given by \eqref{eq:phiK}, for $\ell = 0,\dots,L-1$, regardless of the specific form of $K$.
Let us mention that for some specific choices of $K$ there may be alternative approaches that may lead to an even faster precomputation. Examples include $d$-dimensional rectangles and step hyperbolic crosses, cf. Sections~\ref{Subsec:Rectangle_Precomputation} and \ref{subsec:precomp_stepX}.
As discussed in Section~\ref{se:derivation},
we may write the weights as follows 
\begin{equation*}
        \phi_K(\bsz_\ell) = \frac{1}{N} \sum_{\bsk \in K} \sum_{n=1}^N c_n \omega_{\bsk}(\bsx_n - \bsz_\ell) 
        = \sum_{\bsk \in K} \phih_{\bsk} \omega_{\bsk}(-\bsz_{\ell}),
\end{equation*}
where
\begin{equation*}
    \phih_{\bsk} = \frac{1}{N} \sum_{n=1}^N c_n \omega_{\bsk}(\bsx_n).
\end{equation*}

Let %$M := \abs{\Knugammad}$ and let 
$\bsk_1,\dots,\bsk_{|K|}$ be an enumeration of $K$. We can compute the compression weights above in two steps. For the first step, note that the following matrix-vector product forms an adjoint nonequispaced discrete Fourier transform (NDFT) problem
\begin{equation}\label{eq:ndft}
    \begin{bmatrix}
        \phih_{\bsk_1} \\
        \phih_{\bsk_2} \\
        \vdots \\
        \phih_{\bsk_{|K|}}
    \end{bmatrix} = \frac{1}{N} \begin{bmatrix}
        \omega_{\bsk_1}(\bsx_1) & \omega_{\bsk_1}(\bsx_2) & \hdots & \omega_{\bsk_1}(\bsx_N)\\
        \omega_{\bsk_2}(\bsx_1) & \omega_{\bsk_2}(\bsx_2) & \hdots & \omega_{\bsk_2}(\bsx_N)\\
        \vdots & \vdots & \ddots &\vdots \\
        \omega_{\bsk_{|K|}}(\bsx_1) & \omega_{\bsk_{|K|}}(\bsx_2) & \hdots & \omega_{\bsk_{|K|}}(\bsx_N)\\
    \end{bmatrix} \begin{bmatrix}
        c_1 \\ c_2 \\ \vdots \\ c_N
    \end{bmatrix}.
\end{equation}
Recall that we do not assume that $K$ or $\calX$ have a specific structure. 
If we want to obtain the exact coefficients $\phih_{\bsk}$, then there seems to be no better approach known than the naive one, which is by evaluating the characters $\omega_{\bsk}(\bsx_n)$ and performing the matrix-vector multiplication at cost of order $O(dN \abs{K})$.

In the second step, note that the following matrix-vector product forms the forward discrete Fourier transform problem
\begin{equation}\label{eq:dft}
    \begin{bmatrix}
        \phi_K(\bsz_0) \\ \phi_K(\bsz_1) \\ \vdots \\ \phi_K(\bsz_{L-1})
    \end{bmatrix} = \begin{bmatrix}
        \omega_{\bsk_1}(-\bsz_0) & \omega_{\bsk_2}(-\bsz_0) & \hdots & \omega_{\bsk_{|K|}}(-\bsz_0) \\
        \omega_{\bsk_1}(-\bsz_1) & \omega_{\bsk_2}(-\bsz_1) & \hdots & \omega_{\bsk_{|K|}}(-\bsz_1) \\
        \vdots & \vdots & \ddots & \vdots\\
        \omega_{\bsk_1}(-\bsz_{L-1}) & \omega_{\bsk_2}(-\bsz_{L-1}) & \hdots & \omega_{\bsk_{|K|}}(-\bsz_{L-1}) \\
    \end{bmatrix} \begin{bmatrix}
        \phih_{\bsk_1} \\
        \phih_{\bsk_2} \\
        \vdots \\
        \phih_{\bsk_{|K|}}
    \end{bmatrix}.
\end{equation}
Let $\mathbb{\frakg}$ be the generator of the lattice $\left\{ \bsz_\ell \right\}_{\ell = 0}^{L-1}$. The matrix-vector product above simplifies to a one-dimensional discrete Fourier transform as follows
\begin{equation}\label{univariate_FFT}
    \phi_K(\bsz_\ell) = \sum_{\bsk \in K} \phih_{\bsk} \exp\left( -2\pi i \ell \frac{\bsk \cdot \mathbb{\frakg}}{L}\right) = \sum_{j = 0}^{L-1} \left( \sum_{\substack{ -\bsk \cdot \frakg \equiv j \mod L,\\ \bsk \in K}} \phih_{\bsk} \right) \exp \left( 2\pi i \frac{j\ell}{L} \right).
\end{equation}
First, we compute all values $h_j := \sum_{\bsk \in K, \bsk \cdot \mathbb{\frakg} \equiv -j \mod L} \phih_{\bsk}$ for $j = 0,\dots,L-1$. All the values of $\{h_j\}_{j=0}^{L-1}$ can be computed by iterating once over all elements of $K$. Having computed all $h_j$, the sum above simplifies to
\begin{equation*}
    \begin{aligned}
        \phi_K(\bsz_\ell) &= \sum_{j = 0}^{L-1} h_j \exp \left( 2\pi i \frac{j\ell}{L} \right),
    \end{aligned}
\end{equation*}
which is the one-dimensional inverse fast Fourier transform to evaluate a univariate trigonometric polynomial whose Fourier coefficients are $\{h_j\}_{j=0}^{L-1}$. 
Hence, the cost of solving \eqref{eq:dft} is of the order
\begin{equation*}
  O(  d\abs{K} + L \log(L)),
\end{equation*}
where the hidden constants do not depend on $d$, see \cite{CT1965} and \cite[Algorithm~1]{Kammerer13}.

The algorithm to compute the compression weights using \eqref{eq:ndft} and \eqref{eq:dft} is summarized in \RefAlg{alg:hypX}. 

\SetKwComment{Comment}{/* }{ */}

\begin{algorithm}[H]
    \caption{The data compression algorithm for a general index set $K$}\label{alg:hypX}
    \KwData{$\calX,\calY,\calZ$, K}
    \KwResult{$W_{\calX,\calZ}, W_{\calX,\calY,\calZ}$}
    $(N,d) \gets \textup{size}(\calX), L \gets \textup{size}(\calZ)$ \;
    $A \gets \frac{1}{N}[ \omega_{\bsk}(\bsx) ]_{\bsk \in K, \bsx \in \calX}$ as in \eqref{eq:ndft}\;
    $[\phih_{1,\bsk}]_{\bsk \in K} \gets A \cdot \bsone_N$, $[\phih_{2,\bsk}]_{\bsk \in K} \gets A \cdot \calY$  \tcp*{Expensive matrix-vector products}
    $\left( \bsh_i := [h_{i,j}]_{j=0}^{L-1} \right) \gets \bszero_L $ for $i=1,2$ \;
    \For{$\bsk \in K$}{
        $j = - \bsk \cdot \frakg \mod L$\;
        $h_{i,j} \gets h_{i,j} + \phih_{i,\bsk}$ for $i=1,2$ \;
    }
    $W_{\calX,\calZ} \gets \textsc{IFFT}(\bsh_1) $, $ W_{\calX,\calY,\calZ} \gets \textsc{IFFT}(\bsh_2) $ \tcp*{Fast computations via FFT}
\end{algorithm}

From the discussion above, we see that its total cost is of order
\begin{equation}\label{eq:cost_general} 
    O(d|K|N + L \log(L)).
\end{equation}
This is the cost we have to pay if we deal with the general case 
and cannot exploit structural properties of the index set $K$ or the set
of data points $\calX$.

\begin{remark}\label{Rem:Structure}
As already mentioned, we show in Section~\ref{Subsec:Rectangle_Precomputation} and \ref{subsec:precomp_stepX} how to exploit structure in the cases where $K$ is a $d$-dimensional rectangle or a step hyperbolic cross. 
If, instead, $K$ is a general index set, but  $\calX$ is, e.g., itself a rank-$1$ lattice point set, then we can exploit this additional structure as follows: For the weights $W_{\calX,\calY,\calZ,\ell}$
we may reduce the first step \eqref{eq:ndft} (line 3 of the algorithm above) to a univariate FFT, similarly as in \eqref{univariate_FFT} and solve it at cost $O(d|K| + N \log(N))$, which would be advantageous in the regime where $\log(N)$ is considerably smaller than $d|K|$.
For the weights $W_{\calX, \calZ,\ell}$, i.e., the case where $c_0=\ldots =c_N =1$, the first step is even easier. 
Indeed, if $\mathbb{\frakh}$ denotes the generating vector of $\calX$, then we have 
$\phih_{\bsk} = 1$ if $\bsk \cdot \mathbb{\frakh} \equiv 0 \mod N$ and $\phih_{\bsk} = 0$ otherwise, which can be checked for all $\bsk \in K$ with cost at most $O(d|K|)$.
If, in addition, $\calZ$ is a sub-lattice of $\calX$, i.e., if $\mathbb{\frakg} = j \mathbb{\frakh}$ for some $j\in [N]$, then $\bsk \cdot \mathbb{\frakh} \equiv 0 \mod N$ implies 
$\bsk \cdot \mathbb{\frakg} \equiv 0 \mod N$, leading to $W_{\calX, \calZ,\ell} = |K \cap \calX^{\perp}|$, where $\calX^{\perp}$ is the dual lattice $\{ \bsk \in \Z^d \,:\,\bsk \cdot \mathbb{\frakh} \equiv 0 \mod N\}$ of $\calX$.
Hence, in the sub-lattice case we have for the weights $W_{\calX, \calZ,\ell}$, $\ell = 0, 1, \ldots, L-1$, 
total precomputation cost of order at most $O(d|K| + L)$.
\end{remark}

In Sections \ref{se:contX} to \ref{se:stepX} we consider specific choices of the index set $K$, namely rectangles and two types of hyperbolic crosses.

\subsection{Truncation to a continuous hyperbolic cross}\label{se:contX}

We define the \textbf{weighted continuous hyperbolic cross} for $\nu > 0$, see e.g. \cite{KSW04}, as
\begin{equation}\label{eq:wt_hypX}
    \Knugammad := \left\{ \bsk \in \Z^d : r_\alpha(\bsgamma, \bsk) \leq \nu  \right\}.
\end{equation}
The sets $\Knugammad$ are also referred to in the literature as \textbf{weighted Zaremba crosses}, see e.g. \cite{CKN2010}. %, or as \textbf{symmetric hyperbolic crosses} c.f. \cite{KKP2012,KPV2015}.
For $\bsu \subset [d]$ let 
$$\Knugammau := \left\{ \bsk \in \Z^{d} : \supp \bsk = \bsu \land r_\alpha(\bsgamma,\bsk) \leq \nu \right\},$$
where $\supp \bsk := \{ j\in [d] : k_j \neq 0 \}$.
With this definition, we have the following disjoint union
\begin{equation}\label{hypX_decomp}
    \Knugammad = \bigcup_{\bsu \subseteq [d]} \Knugammau.
\end{equation}
Note that if $\gmu \nu < 1$, then the set $\Knugammau$ is empty.

\subsubsection{Error analysis for weighted Wiener algebras}\label{Sec:cont_Wiener}

\begin{theorem}\label{th:hyp_wtW}
    Let $\alpha > 1, \delta \in (0,\alpha-1), \tau \in (0,\alpha-1-\delta]$ and assume that $g \in \Aalphagammad$. Let $L$ be a prime number, and let $\calZ$ be a lattice point set of size $L$ satisfying \eqref{eq:lattice_wt_err}. Then it holds for $\nu > 1$ that
    \begin{equation*}
        \begin{aligned}
            &\abs{ \frac{1}{N}\sum_{n=1}^N c_n g(\bsx_n) - \frac{1}{L} \sum_{\ell=0}^{L-1} g(\bsz_\ell) \phi_{\Knugammad}(\bsz_\ell) } \leq \err_1(g,\calC) + \err_2(g,\calC)\\
            &\leq \norm{g}_{\Aalphagammad} \left[ \frac{1}{\sqrt{\nu}} + \frac{\sqrt{\nu}}{L^{\alpha-\half-\delta-\tau}} c_{\alpha,\bsgamma,d}\ \zeta_{\delta,d}\ C_{\bsgamma,d}\left(\alpha-\half-\delta,\tau\right)  \right] \mubar_N,
        \end{aligned}
    \end{equation*}
    where $\err_1(g,\calC)$ and $\err_2(g,\calC)$ are defined as in \eqref{eq:error1} and \eqref{eq:error2}, $c_{\alpha,\bsgamma,d}$, $\zeta_{\delta,d}$ and $\mubar_N$ are as in \eqref{eq:const}, and $C_{\bsgamma,d}(\alpha,\tau)$ is as in \eqref{eq:lattice_wt_err}.
\end{theorem}
\begin{proof}
    The first inequality was already established in \eqref{eq:error_gcn}. The bound on $\err_1(g,\calC)$ comes from \eqref{eq:error1} and \RefLem{le:trunc_wtW}, since we have
    \begin{equation*}
        \sup_{\bsk \notin \Knugammad} \frac{1}{\sqrt{r_{\alpha}(\bsgamma, \bsk)}} \leq \frac{1}{\sqrt{\nu}},
    \end{equation*}
    which implies
    \begin{equation*}
        \err_1(g,\calC) \leq \norm{g}_{\Aalphagammad} \frac{\mubar_N}{\sqrt{\nu}}.
    \end{equation*}
    For $\err_2(g,\calC)$, we first note that
    \begin{equation*}
        \max_{\bsk \in \Knugammad} \sqrt{r_{\alpha}(\bsgamma,\bsk)} \leq \sqrt{\nu}.
    \end{equation*}
    Since we have established that $g \phi_{\Knugammad} \in \Halphaminusgamma$ in \RefLem{le:gphi_wtW}, we can use the above relation with \eqref{eq:lattice_wt_err} to get
    \begin{equation*}
        \begin{aligned}
            \err_2(g,\calC) &\leq \norm{g \phi_{\Knud}}_{\Halphaminusgamma} e(\Halphaminusgamma,\calZ) \\
            &\leq \norm{g}_{\Aalphagammad} \sqrt{\nu} c_{\alpha,\bsgamma,d}\ \zeta_{\delta,d} C_{\bsgamma,d}\left(\alpha-\half-\delta,\tau\right) L^{-\alpha+\half+\delta+\tau}\ \mubar_N.
        \end{aligned}
    \end{equation*}
    The statement is hence proved.
\end{proof}

\begin{corollary}\label{cor:wtW_hypX}
    Assume that the conditions of \RefThm{th:hyp_wtW} hold and let $\varepsilon > 0$ be arbitrary. By choosing $\sigma = \sigma(\varepsilon) > 0$ sufficiently small and the parameter $\nu$ as $\nu \asymp L^{(\alpha-\half-\sigma)}$, we get the error bound
    \begin{equation*}
        \err_1(g,\calC) + \err_2(g,\calC) \lesssim L^{-(\alpha - \half)/2+\varepsilon},
    \end{equation*}
    where the implicit constant is independent of $L$.
\end{corollary}

\subsubsection{Error analysis for weighted Korobov spaces}\label{se:wtH}
For the next result, we need the following estimate
\begin{equation}\label{eq:hypX_prop1}
    \sum_{\bsk \in \N^d, \eta(\bsk) > \nu} \eta(\bsk)^{-r} \asymp \nu^{-r+1} (\log(\nu))^{d-1}; \quad r > 1,
\end{equation}
where $\bsk = (k_1,\dots,k_d)$, $\eta(\bsk) := r_{\half}(\bsone_d, \bsk) = \prod_{j \in [d]} \max ( \abs{k_j}, 1 )$, and the implicit constants depend on $r$ and $d$, see e.g. \cite[p.~22]{Dung2018}.

Furthermore, the following estimate on the cardinality of the weighted hyperbolic cross holds for all $\varepsilon > 0$
\begin{equation}\label{eq:hypX_prop2}
    \abs{\Knugammad} \leq \nu^{\frac{1}{2\alpha}+\varepsilon} \widetilde{C}_{\bsgamma,d}(\alpha,\varepsilon), \quad \text{ where } \widetilde{C}_{\bsgamma,d}(\alpha,\varepsilon) := \prod_{j \in [d]} \left( 1 + 2\zeta(1 + 2\alpha \varepsilon ) \gamma_j^{\frac{1}{2\alpha}+\varepsilon} \right),
\end{equation}
see \cite{KSW04}.

\begin{theorem}\label{th:hyp_wtH}
    Let $\alpha > 1, \delta \in (0,\alpha-1), \tau \in (0,\alpha-1-\delta]$ and assume that $g \in \Halphagammad$. Let $L$ be a prime number, and let $\calZ$ be a lattice point set of size $L$ satisfying \eqref{eq:lattice_wt_err}. Then it holds for $\nu > 1$ and $\varepsilon > 0$ that
    \begin{equation*}
        \begin{aligned}
            &\abs{ \frac{1}{N}\sum_{n=1}^N c_n g(\bsx_n) - \frac{1}{L} \sum_{\ell=0}^{L-1} g(\bsz_\ell) \phi_{\Knud}(\bsz_\ell) } \leq \err_1(g,\calC) + \err_2(g,\calC)\\
            &\lesssim \norm{g}_{\Halphagammad} \left[ \frac{\log(\nu)^{\frac{d-1}{2}}}{\nu^{\half-\frac{1}{4\alpha}}} + \frac{\nu^{\half+\frac{1}{4\alpha}+\varepsilon} }{L^{\alpha-\half-\delta-\tau}} c_{\alpha,\bsgamma,d}\ \zeta_{\delta,d} \widetilde{C}_{\bsgamma,d}(\alpha,\varepsilon) C_{\bsgamma,d}\left(\alpha-\half-\delta,\tau\right)  \right] \mubar_N,
        \end{aligned}
    \end{equation*}
    where $\err_1(g,\calC)$ and $\err_2(g,\calC)$ are defined as in \eqref{eq:error1} and \eqref{eq:error2}, $c_{\alpha,\bsgamma,d}$, $\zeta_{\delta,d}$ and $\mubar_N$ are as in \eqref{eq:const}, $C_{\bsgamma,d}(\alpha,\tau)$ as in \eqref{eq:lattice_wt_err}, and $\widetilde{C}_{\bsgamma,d}(\alpha,\varepsilon)$ as in \eqref{eq:hypX_prop2}. Furthermore, the implicit constant in the $\lesssim$-notation depends on $d$, $\alpha$, and $\bsgamma$.
\end{theorem}
\begin{proof}
    The bound on $\err_1(g,\calC)$ comes from \eqref{eq:error1} and by estimating the sum in \RefLem{le:trunc_wtH} as follows. From \eqref{hypX_decomp} we get 
    \begin{equation*}
        \sum_{\bsk \notin \Knugammad} 1/r_{\alpha}(\bsgamma,\bsk) = \sum_{\emptyset \neq \bsu \subseteq [d]} \gmu \sum_{\substack{\bsk \notin \Knugammau,\\ \supp \bsk = \bsu }} \eta(\bsk)^{-2\alpha} =  \sum_{\emptyset \neq \bsu \subseteq [d]} \gmu 2^{\abs{\bsu}} \sum_{\bsk \in \N^{\bsu}, \eta(\bsk)> (\gmu \nu)^{1/2\alpha}} \eta(\bsk)^{-2\alpha}.
    \end{equation*}
    We can now use \eqref{eq:hypX_prop1} to obtain
    \begin{equation*}
        \sum_{\bsk \in \N^{\bsu}, \eta(\bsk)> (\gmu \nu)^{1/2\alpha}} \eta(\bsk)^{-2\alpha} \asymp  (\nu \gmu)^{-1 + 1/2\alpha} 
        %\mig{\left( \frac{1}{2\alpha} \right)^{\abs{\bsu}-1}} 
        \left( \frac{1}{2\alpha} \log (\gmu \nu) \right)^{\abs{\bsu}-1} \leq (\nu \gmu)^{-1 + 1/2\alpha} (\log(\nu))^{\abs{\bsu}-1},
    \end{equation*}
    since $\gmu \leq 1$ and $2\alpha \ge 1$; note that the implicit constants in the $\asymp$-notation depend on $\alpha$ and $d$. Altogether, we have the upper bound
    \begin{equation*}
        \begin{aligned}
            \sum_{\bsk \notin \Knugammad} 1/r_{\alpha}(\bsgamma,\bsk) &\lesssim \frac{1}{\nu^{1 - 1/2\alpha}\log(\nu)} \sum_{\bsu \subseteq [d]} 2^{\abs{\bsu}}  \gmu^{1/2\alpha} (\log(\nu))^{\abs{\bsu}} \\
            &= \frac{1}{\nu^{1 - 1/2\alpha} \log (\nu)} \prod_{j \in [d]} \left( 1 + 2\gamma_j^{1/2\alpha} \log(\nu)  \right) \lesssim \frac{\log(\nu)^{d-1}}{\nu^{1 - 1/2\alpha}},
        \end{aligned}
    \end{equation*}
    and consequently
    \begin{equation*}
        \err_1(g,\calC) \lesssim \norm{g}_{\Halphagammad} \frac{\log(\nu)^{\frac{d-1}{2}}}{\nu^{\half-\frac{1}{4\alpha}}}\ \mubar_N.
    \end{equation*}
    
    For $\err_2(g,\calC)$, we can estimate the sum in \RefLem{le:gphi_wtH} using \eqref{eq:hypX_prop2} as
    \begin{equation*}
        \sum_{\bsk \in \Knugammad} r_{\alpha}(\bsgamma,\bsk) \leq \abs{\Knugammad} \max_{\bsk \in \Knugammad} r_\alpha(\bsgamma,\bsk) \leq \abs{\Knugammad} \nu \leq \\nu^{1+\frac{1}{2\alpha}+ 2\varepsilon} \widetilde{C}_{\bsgamma,d}(\alpha, 2\varepsilon).
    \end{equation*}
    We use the estimate above and \eqref{eq:lattice_wt_err}, as well as \RefLem{le:gphi_wtH}, to get
    \begin{equation*}
        \begin{aligned}
            \err_2(g,\calC) &\leq \norm{g \phiKnugamma}_{\Halphaminusgamma} e(\Halphaminusgamma,\calZ) \\
            &\le \norm{g}_{\Halphagammad} \frac{\nu^{\half+\frac{1}{4\alpha}+\varepsilon} }{L^{\alpha-\half-\delta-\tau}} c_{\alpha,\bsgamma,d}\ \zeta_{\delta,d} \sqrt{ \widetilde{C}_{\bsgamma,d}(\alpha, 2\varepsilon)} C_{\bsgamma,d}\left(\alpha-\half-\delta,\tau\right) \mubar_N.
        \end{aligned}
    \end{equation*}
    The statement is hence proved.
\end{proof}

\begin{corollary}\label{cor:wtH_hypX}
    Assume that the conditions of \RefThm{th:hyp_wtH} hold and let $\varepsilon > 0$ be arbitrary. By choosing $\sigma = \sigma(\varepsilon) > 0$ sufficiently small and the parameter $\nu$ as $\nu \asymp L^{(\alpha-\half-\sigma)}$, we get the error bound
    \begin{equation*}
        \err_1(g,\calC) + \err_2(g,\calC) \lesssim L^{-\left( \frac{\alpha-1}{2}+\frac{1}{8\alpha}\right)+\varepsilon} \log(L)^{\frac{d-1}{2}},
    \end{equation*}
    where the implicit constant is independent of $L$.
\end{corollary}

\begin{remark}\label{re:nu_hypX}
    Note that for both types of spaces, weighted Wiener algebras and weighted Korobov spaces, we chose the parameter $\nu$ of the same order to balance the error. This means that the hyperbolic crosses $\Knugammad$ are of similar size for both spaces.
\end{remark}

\subsubsection{Precomputation of compression weights}\label{se:precomp_wtH}

We now discuss the precomputation of the compression weights with respect to the hyperbolic cross $\Knugammad$ as defined in \eqref{eq:wt_hypX}. 
We rely on the analysis in Section~\ref{se:precomp_gen_K}: There we found that the cost is (at most) of order $O(d |\Knugammad| N + L \log(L))$.
Recall that we are interested in the case where $N$ is (much) larger than $L$. Considering the choice of $\nu \asymp L^{(\alpha-\half-\sigma)}$, see Remark \ref{re:nu_hypX}, and the cardinality of the hyperbolic cross from \eqref{eq:hypX_prop2}, we therefore arrive at a total cost of the order of
\begin{equation*}
  O \left(  dN L^{\half - \frac{1}{4\alpha} + \varepsilon} \widetilde{C}_{\bsgamma,d}(\alpha,\varepsilon) \right),
\end{equation*}
where $\varepsilon > 0$ can be chosen arbitrarily small, and $\widetilde{C}_{\bsgamma,d}(\alpha,\varepsilon) \ge 1$ is as defined in \eqref{eq:hypX_prop2}.
We may improve the precomputation cost to some extent if we are ready to accept a weaker error bound, cf. Section~\ref{Subsec:Comparison_Cost}.

\subsection{Truncation to a high-dimensional rectangle}\label{se:wt_rect}

In this section, we briefly look at another choice of $K$: the high-dimensional rectangle. The main goals are to show in Section~\ref{Subsec:Rectangle_Precomputation} that we may exploit specific structural properties of index sets and to use the findings in the discussion of step hyperbolic crosses in Section~\ref{subsec:precomp_stepX}. 
We define a \textbf{weighted $d$-dimensional rectangle} as
\begin{equation}\label{eq:wt_rect}
    \Rnugammad := \left\{ \bsk \in \Z^d : \max_{j \in [d]} r_{\alpha}(\gamma_j,k_j) \leq \nu \right\}.
\end{equation}

\subsubsection{Error analysis for weighted Wiener algebras}

\begin{theorem}\label{th:wt_rect}
    Let $\alpha > 1, \delta \in (0,\alpha-1), \tau \in (0,\alpha-1-\delta]$ and assume that $g \in \Aalphagammad$. Let $L$ be a prime number, and let $\calZ$ be a lattice point set of size $L$ satisfying \eqref{eq:lattice_wt_err}. Then it holds for $\nu > 1$ that
    \begin{equation*}
        \begin{aligned}
            &\abs{ \frac{1}{N}\sum_{n=1}^N c_n g(\bsx_n) - \frac{1}{L} \sum_{\ell=0}^{L-1} g(\bsz_\ell) \phiRnugamma(\bsz_\ell) } \leq \err_1(g,\calC) + \err_2(g,\calC)\\
            &\leq \norm{g}_{\Aalphagammad} \left[ \frac{1}{\sqrt{\nu}} + \frac{\nu^{d/2}}{L^{\alpha-\half-\delta-\tau}} c_{\alpha,\bsgamma,d}\ \zeta_{\delta,d} C_{\bsgamma,d}\left(\alpha-\half-\delta,\tau\right)  \right] \mubar_N,
        \end{aligned}
    \end{equation*}
    where $\err_1(g,\calC)$ and $\err_2(g,\calC)$ are defined as in \eqref{eq:error1} and \eqref{eq:error2}, $c_{\alpha,\bsgamma,d}$, $\zeta_{\delta,d}$ and $\mubar_N$ are as in \eqref{eq:const}, and $C_{\bsgamma,d}(\alpha,\tau)$ as in \eqref{eq:lattice_wt_err}.
\end{theorem}
\begin{proof}
    We use \eqref{eq:error1} and \RefLem{le:trunc_wtW} to obtain an upper bound $\err_1(g,\calC)$. Hence, we need an estimate on $\sup_{\bsk \notin \Rnugammad} 1/\sqrt{r_{\alpha}(\bsgamma, \bsk)}$.
    If $\bsk \notin \Rnugammad$, then there exists a $j \in [d]$ such that $r_\alpha(\gamma_j,k_j) > \nu$, and for all remaining $j \in [d]$ we have in any case $r_\alpha(\gamma_j,k_j) \geq 1$. Consequently, we get
    \begin{equation*}
        \sup_{\bsk \notin \Rnugammad} \frac{1}{\sqrt{r_{\alpha}(\bsgamma, \bsk)}} \leq \frac{1}{\sqrt{\nu}},
    \end{equation*}
    which implies
    \begin{equation*}
        \err_1(g,\calC) \leq \frac{\norm{g}_{\Aalphagammad}}{\sqrt{\nu}} \mubar_N.
    \end{equation*}
    We know that $\err_2(g,\calC)$ satisfies the inequality
    \begin{equation*}
        \err_2(g,\calC) \leq \norm{g \phiRnugamma}_{\Halphaminusgamma} e(\Halphaminusgamma,\calZ),
    \end{equation*}
    so we may apply \RefLem{le:gphi_wtW} by estimating $\max_{\bsk \in \Rnugammad} \sqrt{r_{\alpha}(\bsgamma,\bsk)}$. Obviously, $\bsk \in \Rnugammad$ satisfies $r_\alpha(\gamma_j,k_j) \leq \nu$ for all $j \in [d]$.
    Therefore, we have
    \begin{equation*}
        \max_{\bsk \in \Rnugammad} \sqrt{r_{\alpha}(\bsgamma,\bsk)} \leq \nu^{d/2}.
    \end{equation*}
    With the estimate above and \eqref{eq:lattice_wt_err}, we finally get
    \begin{equation*}
        \err_2(g,\calC) \leq \norm{g}_{\Aalphagammad} \frac{\nu^{d/2}}{L^{\alpha-\half-\delta-\tau}} c_{\alpha,\bsgamma,d}\ \zeta_{\delta,d}\ C_{\bsgamma,d}\left(\alpha-\half-\delta,\tau\right) \mubar_N.
    \end{equation*}
    The result now follows.
\end{proof}

\begin{corollary}\label{cor:wt_rect}
    Assume that the conditions of \RefThm{th:wt_rect} hold, and let $\varepsilon > 0$ be arbitrary. By choosing  $\sigma = \sigma(\varepsilon)> 0$ sufficiently small and the parameter $\nu$ as $\nu \asymp L^{2(\alpha-\half-\sigma)/(1+d)}$, we get the error bound
    \begin{equation*}
        \err_1(g,\calC) + \err_2(g,\calC) \lesssim L^{-(\alpha-\half)/(1+d)+\varepsilon},
    \end{equation*}
    where the implicit constant is independent of $L$.
\end{corollary}

We see from the above result that as $d$ increases, the error bound gets much worse than the corresponding one for the continuous hyperbolic cross. 
Nevertheless, the advantage in this setting is the fast precomputation of the compression weights, as we will see in Section~\ref{Subsec:Rectangle_Precomputation}.

\subsubsection{Error analysis for weighted Korobov spaces}

For the sake of completeness, we also add error bounds for the Korobov case;  %but we omit the tedious proof details here.
the details of the proof of Theorem~\ref{th:wt_rect_wtH} can be found in the Appendix.

Using \RefLem{le:trunc_wtH} and \RefLem{le:gphi_wtH} for $K = \Rnugammad$, we get the following result.
\begin{theorem}\label{th:wt_rect_wtH}
    Let $\alpha > 1, \delta \in (0,\alpha -1 ), \tau \in (0, \alpha -1 -\delta ]$ and assume that $g \in \Halphagammad$. Let $L$ be a prime number, and let $\calZ$ be a lattice point set of size $L$ satisfying \eqref{eq:lattice_wt_err}. Then it holds for $\nu > 1$ that
    \begin{equation*}
        \err_1(g, \calC) + \err_2(g,\calC) \lesssim \norm{g}_{\Halphagammad}  \left(  \frac{1}{\nu^{\half - \frac{1}{4\alpha}}} + \frac{\nu^{\left(\half + \frac{1}{4\alpha} \right)d}}{L^{\alpha-\half-\delta-\tau}} c_{\alpha,\bsgamma,d}\ \zeta_{\delta,d}\ C_{\bsgamma,d}\left(\alpha-\half-\delta,\tau\right)  \right) 
        \mubar_N,
    \end{equation*}
    where the implicit constant depends on $d$, on $\alpha$, and on the weights $(\gamma_j)_{j\in [d]}$.
\end{theorem}

\begin{corollary}\label{Cor3.14}
    Assume that the conditions of \RefThm{th:wt_rect_wtH} hold and let $\varepsilon > 0$ be arbitrary. By choosing $\sigma = \sigma(\varepsilon) > 0$ sufficiently small and the parameter $\nu$ as $\nu \asymp L^{[(\alpha - \half - \sigma)4\alpha]/[(2\alpha + 1)d + 2\alpha -1]}$, we get the error bound
    \begin{equation*}
        \err_1(g,\calC) + \err_2(g,\calC) \lesssim L^{- \left( [\alpha-\half]/\left[ 1 + \left(\frac{2\alpha+1}{2\alpha - 1}\right) d \right] \right) +\varepsilon},
    \end{equation*}
    where the implicit constant is independent of $L$.
    
\end{corollary}

\subsubsection{Precomputation of the compression weights}\label{Subsec:Rectangle_Precomputation}

The Dirichlet kernel is defined for $n \in \N$ and $x \in [0,1]$ as
\begin{equation}\label{eq:dirichlet}
    D_n(x) := \sum_{k=-n}^{n} \exp (2\pi ikx) = \frac{\sin(2\pi(n+1/2)x)}{\sin(\pi x)}.
\end{equation}
Let $\nu \geq 1$ and $k_j^\ast := \floor{(\gamma_j \nu)^{1/2\alpha}}$ for all $j \in [d]$. Consequently, for $K = \Rnugammad$ as defined in \eqref{eq:wt_rect} we can compute the corresponding compression weights for $\ell = 0,\dots,L-1$, see \eqref{eq:phiK}, efficiently with the help of the ``Dirichlet kernel trick'':
\begin{equation}\label{eq:sin_wts}
    \begin{aligned}[b]
        \phi_{\Rnugammad}(\bsz_\ell) &= \frac{1}{N} \sum_{n=1}^N c_n \sum_{\bsk \in \Rnugammad} \omega_{\bsk}(\bsx_n - \bsz_{\ell}) \\
        &= \frac{1}{N} \sum_{n=1}^N c_n \prod_{j \in [d]} \sum_{k_j = -k_j^\ast}^{k_j^\ast} \exp(2\pi i k_j\cdot(x_{n,j}-z_{\ell,j})) \\
        % &= \frac{1}{N} \sum_{n=1}^N c_n \prod_{j \in [d]} \frac{\sin(2\pi (k_j^\ast + 1/2)(x_{n,j} - z_{\ell,j}) )}{\sin(\pi (x_{n,j}-z_{\ell,j}))}\\
        &= \frac{1}{N} \sum_{n=1}^N c_n \prod_{j \in [d]} D_{k_j^\ast}(x_{n,j} - z_{\ell,j}).
    \end{aligned}
\end{equation}

Hence, the total cost of precomputing all compression weights using the formula above is of the order $O(LNd)$, regardless of the cardinality of the rectangle $\Rnugammad$. 

Alternatively, we may compute the compression weights in two steps: first, the adjoint NDFT problem as in \eqref{eq:ndft}, and second, the DFT problem as in \eqref{eq:dft}. For a high-dimensional rectangle $K = \Rnud$ in the unweighted function space setting, i.e. $\gamma_1 = \dots = \gamma_d = 1$, we know of approximate fast methods for the adjoint NDFT problem that have a computational cost of order $O (\nu^d \log(\nu) + \abs{\log(\overline{\varepsilon})}^dN)$, see \cite{KKP09} and \cite[Algorithm~7.1]{PPST2018}. Here, $\overline{\varepsilon}$ is the desired accuracy of the computation. If we choose $\nu$ to have a balanced error as in \RefCol{cor:wt_rect}, we get a cost of order
\begin{equation*}
   O \left( L^{2(\alpha-\half-\sigma)d/(1+d)} \log(L) + \abs{\log(\overline{\varepsilon})}^dN \right),
\end{equation*}
which grows with $\alpha$, making this approach unfavorable compared to the direct computation via the ``Dirichlet kernel trick'' in \eqref{eq:sin_wts}.

\subsection{Truncation to a step hyperbolic cross}\label{se:stepX}

We define a \textbf{weighted step hyperbolic cross} as
\begin{equation}\label{eq:stepX}
    \Qmgammad := \bigcup_{\substack{ \bst \in \N_{0}^{d},\\ \norm{\bst}_1 = m }} \left\{ \bsk \in \Z^d : r_{\alpha}(\gamma_j, k_j) \leq 2^{t_j} \quad \forall j \in [d] \right\},
\end{equation}
see e.g. \cite[Section~2.3]{Dung2018}. The sets $\Qmgammad$ are also known as \textbf{dyadic hyperbolic crosses} in the literature; see, e.g., \cite{KKP2012}.
For $\bss \in [1,\infty)^d$, let
\begin{equation*}
    I^{\alpha}_{\bss, \bsgamma,d} := \{ \bsk \in \Z^d : r_\alpha(\gamma_j, k_j) \leq s_j \quad \forall j \in [d] \}.
\end{equation*}
Using the intervals above, we have the following equivalent definition for the step hyperbolic cross
\begin{equation}\label{eq:stepX_equiv1}
    \Qmgammad = \bigcup_{\substack{ \bss \in \N^d,\\ \prod_{j \in [d]} s_j = 2^m }} I^{\alpha}_{\bss, \bsgamma,d},
\end{equation}
due to the fact that $\prod_{j \in [d]} s_j = 2^m$ for $m \in \N_0$ and $\bss \in \N^d$ implies, via the uniqueness of the prime factorization, that for all $j \in [d]$ there exists some $t_j \in \N_0$ such that $s_j = 2^{t_j}$.

\begin{lemma}\label{le:stepX}
    Let $\alpha > 1/2$ and $m \in \N_0$. Then we have
    \begin{equation*}
        K^{\alpha}_{2^{m-(d-1)},\bsgamma,d} \subseteq \Qmgammad \subseteq K^{\alpha}_{2^{m},\bsgamma,d}.
    \end{equation*}
\end{lemma}

Note that for $m < d-1$ the set $K^{\alpha}_{2^{m-(d-1)},\bsgamma,d}$ is actually empty, i.e., the first inclusion in Lemma~\ref{le:stepX} is trivially satisfied. Although we are only interested in the case $\alpha > 1/2$, the proof below reveals  that the statement of the lemma is still true if we allow arbitrary $\alpha >0$.

\begin{proof}
    We first show $\Qmgammad \subseteq K^{\alpha}_{2^{m},\bsgamma,d}$. We claim that the continuous hyperbolic cross, as defined in \eqref{eq:wt_hypX}, has the following equivalent representation
    \begin{equation}\label{eq:wt_hypX_equiv}
        K^{\alpha}_{2^{m},\bsgamma,d} = \bigcup_{\substack{ \bss \in [1,\infty)^d,\\ \prod_{j \in [d]} s_j = 2^m }} I^{\alpha}_{\bss, \bsgamma,d}.
    \end{equation}
    Let us now verify this claim. It is clear that
    \begin{equation*}
        \left\{ \bsk \in \Z^d : \prod_{j \in [d]} r_{\alpha}(\gamma_j, k_j) \leq 2^m \right\} \supseteq  \bigcup_{\substack{ \bss \in [1,\infty)^d,\\ \prod_{j \in [d]} s_j = 2^m }} \{ \bsk \in \Z^d : r_{\alpha}(\gamma_j, k_j) \leq s_j \quad \forall j \in [d] \}.
    \end{equation*}
    In order to show the reverse inclusion, let $\bsk \in K^{\alpha}_{2^m,\bsgamma,d}$. Let $r_j := r_{\alpha}(\gamma_j, k_j) \in [1,\infty)$ for all $j \in [d]$. Now put $s_j := r_j$ for $j = 1,\dots,d-1$, and $s_d := 2^m / \prod_{j \in [d-1]} r_j$. Since $\prod_{j \in [d]} r_j = r_\alpha(\bsgamma, \bsk) \leq 2^m$, we have $r_j \leq s_j$ for all $j \in [d]$, and $\prod_{j \in [d]} s_j = 2^m$. Hence, $\bsk \in I^{\alpha}_{\bss, \bsgamma,d}$ for $\bss = \{ s_j \}_{j \in [d]}$, and we have shown \eqref{eq:wt_hypX_equiv}. The inclusion $\Qmgammad \subseteq K^{\alpha}_{2^{m},\bsgamma,d}$ now follows from \eqref{eq:stepX_equiv1} and \eqref{eq:wt_hypX_equiv}.

    We now show $K^{\alpha}_{2^{m-(d-1)},\bsgamma,d} \subseteq \Qmgammad$ via induction on $d$. For $d=1$, we have 
    \begin{equation*}
        K^{\alpha}_{2^{m},\gamma_1,1} = \{ k \in \Z : r_{\alpha}(\gamma_1, k) \leq 2^{m} \} = I^{\alpha}_{2^{m},\gamma_1,1} = Q^{\alpha}_{m,\gamma_1,1}.
    \end{equation*}
    Now, let $d \geq 2$ and write $\bsgamma^{\prime} := \{ \gamma_j \}_{j \in [d-1]}$. We assume for the induction hypothesis that the relation $K^{\alpha}_{2^{m^{\prime}-(d-2)},\bsgamma^{\prime},d-1} \subseteq Q^{\alpha}_{m^{\prime},\bsgamma^{\prime},d-1}$ holds for all  $m' \in \N_0$. 
    If $m < d-1$, then $K^{\alpha}_{2^{m-(d-1)},\bsgamma,d} =  \emptyset$ and there remains nothing to show. So let $m \ge d-1$ and $\bsk \in K^{\alpha}_{2^{m-(d-1)},\bsgamma,d} \neq \emptyset$. Then, due to \eqref{eq:wt_hypX_equiv}, there exists an $\bss \in [1,\infty)^{d}$ such that $\bsk \in I^{\alpha}_{\bss, \bsgamma,d}$ and $\prod_{j \in [d]} s_j = 2^{m-(d-1)}$. Now put $\sigma_d := \log_2(s_d)$ and $t_d := \ceil{\sigma_d}$, which implies $\sigma_d \geq t_d-1$. Then, we have
    \begin{equation*}
        \prod_{j \in [d-1]} s_j = 2^{m - \sigma_d - (d-1)} \leq 2^{m - (t_d-1) - (d-1)},
    \end{equation*}
    which implies, using the induction hypothesis, for $\bsk^{\prime} := \{ k_j \}_{j \in [d-1]}$
    \begin{equation*}
        \begin{aligned}
            \bsk^\prime \in K^{\alpha}_{2^{m-t_d-d+2},\bsgamma^{\prime},d-1} &\subseteq Q^{\alpha}_{m-t_d,\bsgamma^{\prime},d-1}.
        \end{aligned}
    \end{equation*}
    Note that we may rewrite the step-hyperbolic cross as
    \begin{equation*}
        \Qmgammad = \bigcup_{\tau_d = 0}^{m} \left\{ \bsell \in \Z^d : r_{\alpha}(\gamma_d,\ell_d) \leq 2^{\tau_d} \land \bsell^\prime \in Q^{\alpha}_{m-\tau_d,\bsgamma^{\prime},d-1} \right\}.
    \end{equation*}
    The last two formulas, and the inequality $r_\alpha(\gamma_d,k_d) \leq s_d = 2^{\sigma_d} \leq 2^{t_d}$ are now sufficient to conclude that $\bsk \in \Qmgammad$. Hence, the result follows.

\end{proof}

\begin{figure}
    \centering
    \includegraphics[width=400pt]{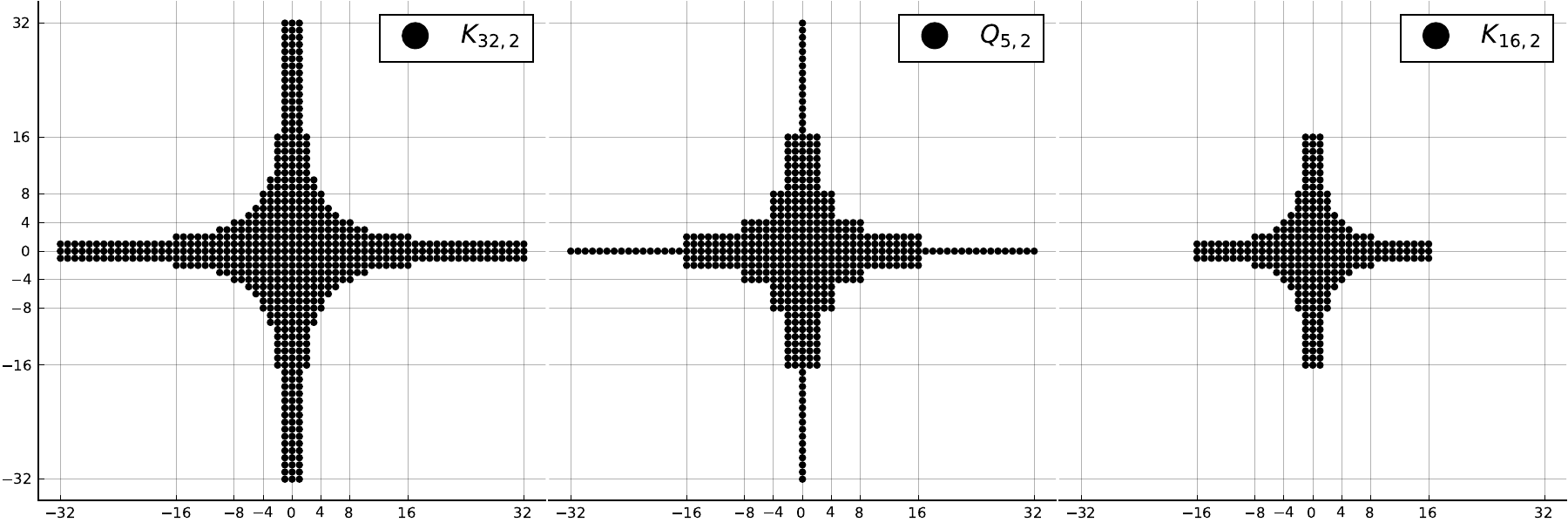}
    \caption{The hyperbolic crosses $K^{\alpha}_{2^{m},\bsgamma,d}$, $\Qmgammad$ and  $K^{\alpha}_{2^{m-(d-1)},\bsgamma,d}$ for $m=5$, $\bsgamma = (1,1)$ and $d=2$, labeled as $K_{32,2}$, $Q_{5,2}$ and $K_{16,2}$ respectively.}
    \label{fig:hyperbolic_cross}
\end{figure}

\begin{remark}
    In general, we do not have $\Qmgammad = K^{\alpha}_{2^{m},\bsgamma,d}$, and the step hyperbolic cross is a proper subset of the continuous hyperbolic cross. We illustrate the hyperbolic crosses for $\alpha = 0.5, d=2, \bsgamma = (1,1)$ and $m=5$ in \RefFig{fig:hyperbolic_cross}. Take, for example, $\bsk = (6,5)$. Then $r_\alpha(\bsgamma,\bsk) = 30 \leq 2^5$, which means $\bsk \in K^{\alpha}_{2^{5},\bsgamma,2}$. However, $\bsk \notin Q^\alpha_{5,\bsgamma,2}$ since $2^2 < r_\alpha(1,6),r_\alpha(1,5) \leq 2^3$, i.e. $\bsk \in I^{\alpha}_{(8,8),\bsgamma} \setminus Q^\alpha_{5,\bsgamma,2}$.
\end{remark}

\begin{remark}
    Similar inclusions as in \RefLem{le:stepX} can be found in the literature; see, e.g., \cite[Lemma~2.6]{KPV2015} and \cite[Lemma~2.1]{KKP2012}.
\end{remark}

\subsubsection{Error analysis for weighted Wiener algebras}\label{Sec:step_Wiener}

We may use \RefLem{le:stepX} and reduce the error analysis in the case where we choose our index sets as step hyperbolic crosses
to the one in Section~\ref{Sec:cont_Wiener}, where we choose them as continuous hyperbolic crosses.

\begin{theorem}\label{th:stepX}
    Let $\alpha > 1, \delta \in (0,\alpha-1), \tau \in (0,\alpha-1-\delta]$ and assume that $g \in \Aalphagammad$. Let $L$ be a prime number, and let $\calZ$ be a lattice point set of size $L$ satisfying \eqref{eq:lattice_wt_err}. Then it holds for $m \in \N$ that
    \begin{equation*}
        \begin{aligned}
            &\abs{ \frac{1}{N}\sum_{n=1}^N c_n g(\bsx_n) - \frac{1}{L} \sum_{\ell=0}^{L-1} g(\bsz_\ell) \phiQmgamma(\bsz_\ell) } \leq \err_1(g,\calC) + \err_2(g,\calC)\\
            &\leq \norm{g}_{\Aalphagammad} \left[ \frac{1}{\sqrt{2^{m-d+1}}} + \frac{\sqrt{2^{m}}}{L^{\alpha-\half-\delta-\tau}} c_{\alpha,\bsgamma,d}\ \zeta_{\delta,d}\ C_{\bsgamma,d}\left(\alpha-\half-\delta,\tau\right)  \right]  \mubar_N,
        \end{aligned}
    \end{equation*}
    where $\err_1(g,\calC)$ and $\err_2(g,\calC)$ are defined as in \eqref{eq:error1} and \eqref{eq:error2}, 
    $c_{\alpha,\bsgamma,d}$, $\zeta_{\delta,d}$ and $\mubar_N$ are as in \eqref{eq:const}, and $C_{\bsgamma,d}(\alpha,\tau)$ as in \eqref{eq:lattice_wt_err}.
\end{theorem}
\begin{proof}
    We use again \eqref{eq:error1} and \RefLem{le:trunc_wtW} to bound $\err_1(g,\calC)$. Hence, we need an estimate on $\sup_{\bsk \notin \Qmgammad} 1/\sqrt{r_{\alpha}(\bsgamma, \bsk)}$. Since $\Qmgammad \supseteq K^\alpha_{2^{m-d+1},\bsgamma,d}$, see \RefLem{le:stepX}, we get
    \begin{equation*}
        \sup_{\bsk \notin \Qmgammad} \frac{1}{\sqrt{r_{\alpha}(\bsgamma, \bsk)}} \leq \sup_{\bsk \notin K^{\alpha}_{2^{m-d+1}, \bsgamma, d}} \frac{1}{\sqrt{r_{\alpha}(\bsgamma, \bsk)}} \leq \frac{1}{\sqrt{2^{m-d+1}}},
    \end{equation*}
    which means
    \begin{equation*}
        \err_1(g,\calC) \leq \frac{\norm{g}_{\Aalphagammad}}{\sqrt{2^{m-d+1}}}\ \mubar_N.
    \end{equation*}
    For $\err_2(g,\calC)$, we may apply \RefLem{le:gphi_wtW} by estimating $\max_{\bsk \in \Qmgammad} \sqrt{r_{\alpha}(\bsgamma,\bsk)}$. Since $\Qmgammad \subseteq K^{\alpha}_{2^{m},\bsgamma,d}$, see \RefLem{le:stepX}, we get
    \begin{equation*}
        \max_{\bsk \in \Qmgammad} \sqrt{r_{\alpha}(\bsgamma,\bsk)} \leq \max_{\bsk \in K^{\alpha}_{2^{m},\bsgamma,d}} \sqrt{r_{\alpha}(\bsgamma,\bsk)} \leq \sqrt{2^{m}}.
    \end{equation*}
    With the estimate above and \eqref{eq:lattice_wt_err}, we finally get
    \begin{equation*}
        \begin{aligned}
            \err_2(g,\calC) &\leq \norm{g \phiQmgamma}_{\Halphaminusgamma} e(\Halphaminusgamma,\calZ) \\
            \\ &\leq \norm{g}_{\Aalphagammad} \frac{\sqrt{2^{m}}}{L^{\alpha-\half-\delta-\tau}} c_{\alpha,\bsgamma,d}\ \zeta_{\delta,d}\ C_{\bsgamma,d}\left(\alpha-\half-\delta,\tau\right) \mubar_N.
        \end{aligned}
    \end{equation*}
    The result now follows.
\end{proof}

\begin{corollary}\label{cor:wtW_stepX}
    Assume that the conditions of \RefThm{th:stepX} hold, and let  $\varepsilon > 0$ be arbitrary. By choosing $\sigma = \sigma(\varepsilon) > 0$ sufficiently small and the parameter $m$ as $m = \left(\alpha-\half-\sigma \right) \log_2(L) + (d-1)/2$, we get the error bound
    \begin{equation}\label{eq:stepX_err}
        \err_1(g,\calC) + \err_2(g,\calC) \lesssim L^{-\left(\alpha - \half\right)/2 + \varepsilon},
    \end{equation}
    where the implicit constant is independent of $L$.
\end{corollary}

\subsubsection{Error analysis for weighted Korobov spaces}

We have shown in \RefLem{le:stepX} how the step hyperbolic cross is related to the continuous one. 
Similarly as in the case of weighted Wiener algebras in Section~\ref{Sec:step_Wiener}, 
we may use \RefLem{le:stepX} in the case of weighted Korobov spaces and play the error analysis for step hyperbolic crosses
back to the one for hyperbolic crosses. 
Therefore, we may conclude that if we use the truncation to a step hyperbolic cross $\Qmgammad$, with $g \in \Halphagammad$, then, we may choose $m$ as in \RefCol{cor:wtW_stepX} to get the same error bound as mentioned in \RefCol{cor:wtH_hypX}.

\subsubsection{Precomputation of compression weights}\label{subsec:precomp_stepX}
To facilitate a fast computation of the compression weights, we represent the step hyperbolic cross as a suitable disjoint union of rectangles
\begin{equation}\label{eq:stepX_union}
    \Qmgammad = \bigcup_{\substack{ \bst \in \N_0^d \\ \norm{\bst}_1 =  m }} \left\{ \bsk \in \Z^d : r_\alpha(\gamma_1, k_1) \leq 2^{t_1} \, \wedge \, \floor{2^{t_j - 1}} < r_\alpha(\gamma_j, k_j) \leq 2^{t_j} \quad \text{for} \quad j =2, \ldots, d \right\}.
\end{equation}
When deriving the formula for $\phiQmgamma(\bsz_\ell)$, it will become apparent that treating the first coordinate differently from the others allows us to reduce the computational cost of the weights by reducing the number of factors under the product sign by one.
In the representation \eqref{eq:stepX_union}, we call the vectors $\bst$ shape vectors. For $j \in [d]$ and $t_j \in \N_0$, define
\begin{equation*}
    \taulow_j := \floor{ \left( \gamma_j 2^{t_j-1} \right)^{1/2\alpha} } \quad \text{and} \quad \tauup_j := \floor{ \left( \gamma_j 2^{t_j} \right)^{1/2\alpha} }.
\end{equation*}
Note that for $t_j = 0$, we have $\taulow_j = 0$, as well as $\tauup_j = 0$ for $\gamma_j < 1$, and $\tauup_j = 1$ for $\gamma_j = 1$. We can now express the compression weights as follows
\begin{equation}\label{eq:stepX_weights}
    \begin{aligned}
        \phiQmgamma(\bsz_\ell) = &\frac{1}{N} \sum_{n=1}^N c_n \sum_{\bsk \in \Qmgammad} \omega_{\bsk}(\bsx_n - \bsz_{\ell}) \\
        = &\frac{1}{N} \sum_{n=1}^N c_n \sum_{\substack{\bst \in \N_{0}^{d},\\ \norm{\bst}_1 =  m}} \left( \sum_{ \substack{k_1 \in \Z :\\ r_\alpha(\gamma_1,k_1) \leq 2^{t_1}}}\exp(2\pi i k_1 (x_{n,1}-z_{\ell,1})) \right) \times\\
        &\times \prod_{j=2}^d \sum_{ \substack{k_j \in \Z:\\ \floor{2^{t_j-1}} < r_\alpha(\gamma_j,k_j) \leq 2^{t_j}}} \exp(2\pi i k_j (x_{n,j}-z_{\ell,j})) \\
        = &\frac{1}{N} \sum_{n=1}^N c_n \sum_{\substack{\bst \in \N_{0}^{d},\\ \norm{\bst}_1 =  m}} \left( \sum_{k_1 = -\tauup_1}^{\tauup_1} \exp(2\pi i k_1 (x_{n,1}-z_{\ell,1})) \right)\times\\
        &\times \prod_{j=2}^d  \left( \sum_{k_j = -\tauup_j}^{\tauup_j} \exp(2\pi i k_j (x_{n,j}-z_{\ell,j}))  - \sum_{k_j = -\taulow_j}^{\taulow_j} \exp(2\pi i k_j (x_{n,j}-z_{\ell,j})) \right) \\
        = &\frac{1}{N} \sum_{n=1}^N c_n \sum_{\substack{\bst \in \N_{0}^{d},\\ \norm{\bst}_1 =  m
          }} D_{\tauup_1}(x_{n,1}-z_{\ell,1}) \prod_{j =2}^d \left( D_{\tauup_j}(x_{n,j}-z_{\ell,j}) - D_{\taulow_j}(x_{n,j}-z_{\ell,j}) \right) ,
    \end{aligned}
\end{equation}
where the last line follows from the definition of the Dirichlet kernel in \eqref{eq:dirichlet}.

The algorithm that uses the last expression in \eqref{eq:stepX_weights} to calculate the compression weights $\phiQmgamma(\bsz_\ell)$
has the following total computing time:
For a given $\bst$, the product appearing in the final expression of \eqref{eq:stepX_weights}, whose factors are a single Dirichlet kernel and the differences of Dirichlet kernels, can be computed in $O(d)$ operations.
The set of shape vectors 
\begin{equation}
    T(m,d) := \left\{ \bst \in \N_{0}^{d} : \norm{\bst}_1 =  m \right\}
\end{equation}
has cardinality 
\begin{equation*}
    |T(m,d)| =  \binom{d-1+m}{d-1} \leq \min \left( \frac{(d-1+m)^m}{m!}, \frac{(d-1+m)^{d-1}}{(d-1)!} \right).
\end{equation*}
Hence, the total cost of computing the compression weights $\phiQmgamma(\bsz_\ell)$ for $\ell = 0,\dots,L-1$ is of the order of
\begin{equation}\label{precomp_time:stepX}
  O \left(  \min \left( \frac{(d-1+m)^m}{m!}, \frac{(d-1+m)^{d-1}}{(d-1)!} \right) NLd \right).
\end{equation}

\begin{figure}
    \centering
    \includegraphics[width=400pt]{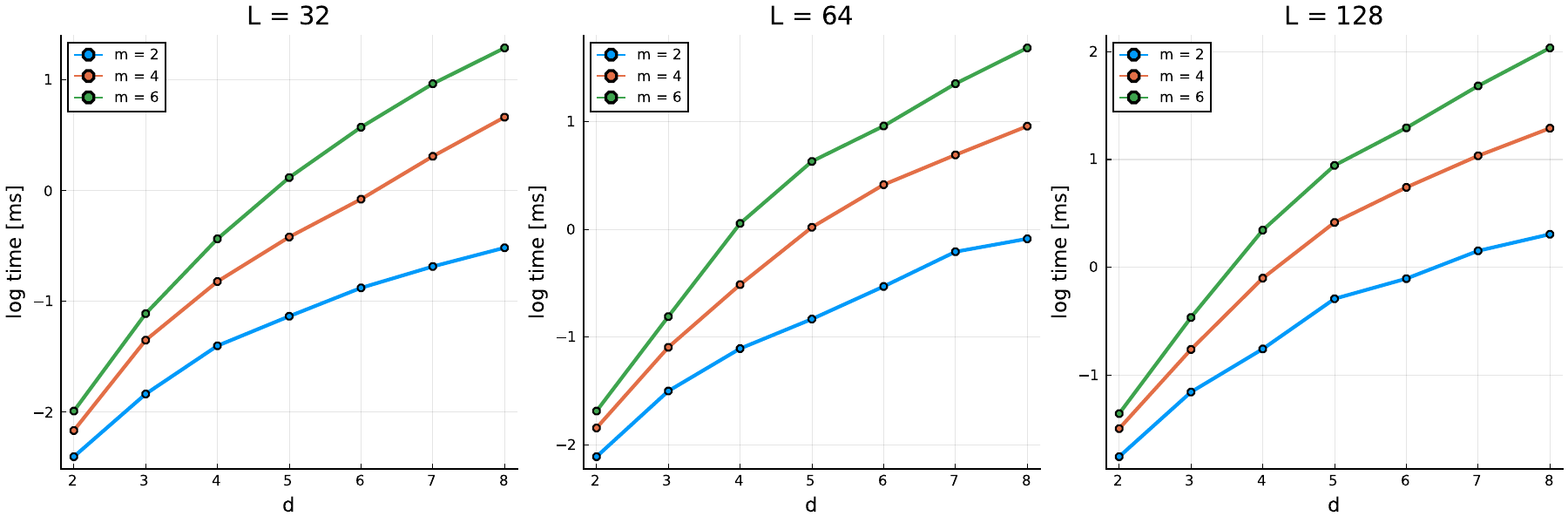}
    \caption{Benchmark tests for $N = 1000$ samples generated uniformly random on $[0,1]^d$ with responses in $[0,1]$ using the step hyperbolic cross.}
    \label{fig:times_stepX}
\end{figure}

To give an impression of the actual computing times of the algorithm based on \eqref{eq:stepX_weights}, 
we prepared some
small numerical examples, see Figure~\ref{fig:times_stepX}.
For dimensions $d=2,\ldots,8$, we generated $N=1000$ i.i.d. random samples from the uniform distribution on $[0,1]^d$ (one sample for each
dimension) and corresponding i.i.d. responses from the uniform distribution on $[0,1]$. Then we  compressed the resulting data sets $\calX$ to lattice point sets $\calZ$ consisting of $32$, $64$, and $128$ points, respectively. 
The lattice point sets are rank-$1$ lattices, whose generating vectors were also chosen randomly by sampling without replacement from the set of coprimes of $L$.
We calculated the corresponding compression weights for $\Qmgammad$ for $m = 2,4,6$, %$\Knugammad$, 
choosing $\alpha = 1.001$ and $\bsgamma = \bsone_d$. To obtain a fast algorithm using broadcasting for matrix operations, we can arrange \eqref{eq:stepX_weights} into a matrix-matrix product which is fast with linear algebra libraries. The values of $\{\tauup_j\}_{j \in [d]}$ and $\{\taulow_j\}_{j \in [d]}$ for all shape vectors can be stored for a slight improvement.
The precomputation times listed in Figure~\ref{fig:times_stepX} are the time it
took on an 11th Gen Intel\textsuperscript{\textregistered} Core\texttrademark i5-1140G7 processor with 16GB RAM
to generate the whole sets of weights $\WXZl$ and $\WXYZl$, $\ell =1, \ldots, L$.
These results suggest that the precomputation time of the compression weights grows rather sub-exponentially in the dimension $d$ than exponentially.
%as indicated  by the theoretical estimate \eqref{precomp_time:stepX}.

Let us close this subsection with the short remark that an alternative to the algorithm described above is to use the two-step procedure as in \RefAlg{alg:hypX}, where we could use the fast approximate Fourier transform for the first step, mentioned in \eqref{eq:ndft}. This transform for the step hyperbolic cross is known as the nonequispaced sparse FFT (NSFFT), and depends on the representation involving the union over high-dimensional rectangles \eqref{eq:stepX_equiv1}, see \cite[Section~4.2]{KKP09} and \cite{FKP2006}. However, the total cost would still be unfavorable, see the discussion after \eqref{eq:sin_wts}. Moreover, the NSFFT implementation is available only for $d=2,3$.

\section{Comparison to data compression with digital nets}\label{se:comparison}

In this section, we want to compare our compression algorithms based on lattice point sets with the ones of Dick and Feischl that are based on digital $(t,m,d)$-nets, see \cite{DF2020}. 
In their error analysis of higher-order digital nets, see \cite[Sect.~4.2.2]{DF2020}, they consider as function spaces the Sobolev spaces $W_{2,\alpha,d}$, where the order $\alpha \ge 2$ is a natural number.
More precisely,
the multivariate Hilbert space $W_{2,\alpha,d}$, $\alpha \in \N$, is the $d$-fold tensor product of the univariate unanchored Sobolev space defined as
\begin{equation*}
    W_{2,\alpha,1} = \left\{ g \in L^2([0,1]) : \forall j \in [\alpha] : g^{(j)} \in L^2([0,1]) \right\},
\end{equation*}
endowed with the norm
\begin{equation*}
    \norm{g}_{W_{2,\alpha,1}}^2 = \sum_{j = 0}^{\alpha - 1} \abs{\int_{[0,1]}  g^{(j)}(x) \dx}^2 + \int_{[0,1]} \abs{ g^{(\alpha)}(x)}^2 \dx.
\end{equation*}
Note that for $g\in H_{\alpha, \gamma, 1}$ the summands for $j=1, \ldots, \alpha -1$ vanish due to (Lebesgue's version of) the fundamental theorem of calculus and the periodicity of $g, g', \ldots, g^{(\alpha-2)}$.
Hence, due to \eqref{korobov_norm_derivatives}, the Korobov space $H_{\alpha, \gamma, 1}$ is continuously embedded into $W_{2,\alpha,1}$, and for the choice $\gamma = 1/(2\pi)^{2\alpha}$ 
the embedding constant (i.e., the operator norm of the identical embedding)  is given by $1$. The same holds for the function spaces of $d$-variate functions $\Halphagammad$ and $W_{2,\alpha,d}$ if we choose
\begin{equation}\label{gamma_2_pi_alpha}
    \gamma_1 = \gamma_2 = \dots = \gamma_d = (2\pi)^{-2\alpha},
\end{equation}
since these spaces are the $d$-fold Hilbert space tensor products of $H_{\alpha, \gamma, 1}$ and $W_{2,\alpha,1}$, respectively. We will stick to this choice of (equal) weights and denote the corresponding Korobov space $\Halphagammad$ (admittedly, with some abuse of notation) by $\Halphad$, and suppress the reference to the weights $\bsgamma$ also in other notation.
Although $W_{2,\alpha,d}$ is  larger than $\Halphad$, since it contains in addition certain non-periodic functions, we believe that these choices of function spaces allow essentially for a fair comparison of the data compression approaches based on lattice point sets
and digital $(t,m,d)$-nets. Indeed, in many results for numerical integration and approximation, the convergence rates on both types of spaces (``periodic and non-periodic Sobolev spaces'') are the same.

For the analysis of plain digital $(t,m,d)$-nets, Dick and Feischl considered Sobolev spaces $W_{1,2,d}$
of functions with dominating mixed smoothness of order $\alpha =2$, endowed with
the norm 
\begin{equation*}
    \norm{g}_{W_{1,2,d}} := \sum_{\bsu \subseteq [d]} \sum_{\bsv \subseteq \bsu} %\sum_{\bstau \in [\alpha-1]^{\abs{\bsu \setminus \bsv}}} 
    \int_{[0,1]^{\bsv}} \abs{ \int_{[0,1]^{d-\abs{\bsv}}} \left( \prod_{j \in \bsv} \partial_{z_j}^{2} \prod_{j \in \bsu \setminus \bsv} \partial_{z_j} \right) g(\bsz) \rd \bsz_{[d]\setminus \bsv} } \rd \bsz_\bsv,
\end{equation*}
where the inner integral is over all variables $z_i$, $i\in [d]\setminus \bsv$,
and in the outer integral over all $z_j$, $j\in \bsv$, see \cite[Sect.~4.2.1]{DF2020}. Since these spaces are larger than the ones we considered previously in this paper, they do not allow for a fair comparison. Indeed, our results for Wiener algebras or Korobov spaces provide a substantially higher convergence rate than the ones presented in \cite[Thm.~11 \& Cor.~12]{DF2020}. 
That is why we stick to the Hilbert space case and compare the results from \cite[Sect.~4.2.2]{DF2020} with our results for Korobov spaces $\Halphad$.

\subsection{Comparison of errors}

In order to have a similar notation for the comparison of lattices with digital $(t,m,d)$-nets in base $b$, we set $L = b^m$. This implies $m \asymp \log(L)$.
Assume that $g \in W_{2,\alpha,d}$ for some integer $\alpha \geq 2$. Let $\calZ$ be an order $\alpha$ digital $(t, m, d)$-net in base $b$. Then, as shown in \autocite[Theorem~13]{DF2020}, it holds for $m - t \geq \nu \in \N$
\begin{equation*}
    \err_1(g,\calC) + \err_2(g,\calC) \lesssim \norm{g}_{W_{2,\alpha,d}} \left( \nu^{2d-1} b^{-\nu} + m^{\alpha d} b^{-\alpha(m-\nu)} \right) \frac{1}{N}\sum_{n=1}^{N} \abs{c_n}.
\end{equation*}
(Note that the previous inequality differs slightly from \cite[Theorem~13]{DF2020}, 
since Dick and Feischl did not take into account the whole error term $\err_1(g,\calC)$ as defined in~\eqref{eq:error1},
but only the term $\| g-g_K \|_{\infty}$, instead, cf. \cite[Formula~(17)]{DF2020}.)
One may now choose $\nu = \lceil \frac{\alpha}{\alpha+1}m \rceil$ to balance the powers of $b$ in the error terms above to get
\begin{equation}\label{eq:nets_p2}
    \err_1(g,\calC) + \err_2(g,\calC) \lesssim m^{\alpha d} b^{-\alpha m / (\alpha + 1)} \asymp L^{-\frac{\alpha}{\alpha+1}} \log(L)^{\alpha d}, 
\end{equation}
as stated in \cite[Remark~15]{DF2020}.

Let $\alpha > 1$ %, $\delta \in (0,\alpha-1)$ 
and assume that $g \in \Halphad$. Let $\calZ$ be a lattice point set of size $L$, where $L$ is a prime number. 
We denote the continuous hyperbolic cross based on equal weights as in \eqref{gamma_2_pi_alpha} by $\Knud$. For an $\varepsilon> 0$ that can be chosen arbitrarily small, we may select the parameters of our compression algorithm based on $L$ and $\Knud$ as in Corollary~\ref{cor:wtH_hypX} 
to get the error bound
\begin{equation}\label{eq:lat_p2}
    \err_1(g,\calC) + \err_2(g,\calC) \lesssim L^{-\left(\frac{\alpha - 1}{2}+\frac{1}{8\alpha}\right)+\varepsilon} \log (L)^{\frac{d-1}{2}},
\end{equation}
where the implicit constant is independent of $L$.

We conclude that the compression algorithm based on lattices and continuous hyperbolic crosses offers a better bound on the convergence rate than the one based on digital nets for all integers $\alpha > 2$.
Note, in particular, that the guaranteed (polynomial) convergence rate of the error in \eqref{eq:nets_p2} is $\frac{\alpha}{\alpha+1}$, which is strictly smaller than $1$ for all values of $\alpha$, while the convergence rate in \eqref{eq:lat_p2} is strictly larger than $1$ for any integer $\alpha > 2$, increasing essentially linearly in $\alpha$,  and can be arbitrarily high for sufficiently smooth functions.

\subsection{Comparison of precomputation cost}\label{Subsec:Comparison_Cost}

We now compare the cost of computing the compression weights for the two approaches.
We again emphasize that this cost belongs to the ``realm of precomputation'', since these weights can be used over and over again to evaluate the approximate error \eqref{eq:error_app} for arbitrary parameters $\theta$. Recall from the previous section, that the optimal choice of $\nu \in \N$ for digital $(t,m,d)$-nets is 
$\ceil{m\alpha/(\alpha+1)}$ for $\alpha \geq 2$.
The precomputation cost of evaluating the compression weights $\WXZl,\WXYZl$ for $\ell = 0,\dots,b^m-1$ using the algorithm described in \autocite*[Algo~2,Algo~3,Lemma~4,Algo~6,Lemma~7]{DF2020} is of the order of
\begin{equation}\label{DF:first_cost_bound}
    d^2 m b^m N \asymp d^2 NL \log(L),
\end{equation}
see also \cite[Theorem~8]{DF2020}. A modified version of that algorithm using digital nets for $\alpha \geq 2$  was described in \cite[Lemma~5]{DF2020}, in case we know an upper bound on the value $t$ for the digital net. In this approach, the number of operations is reduced by using additional storage of size $N$.
Then the cost of computing the compression weights with this modified approach for $\nu = \ceil{m\alpha/(\alpha+1)}$ and $t \leq m/(\alpha+1)$ is of the order
\begin{equation}\label{DF:second_cost_bound}
    \binom{d-1+\nu}{d-1} N b^{m/(\alpha+1)+d-1} \asymp \binom{d-1+\nu}{d-1} b^{d-1} NL^{1/(\alpha+1)},
\end{equation}
where $m = \log_b(L) \in \N$, see \cite[Lemma~5]{DF2020}.
Let us again stress that we assume $L$ to be relatively small compared to $N$. That is why
the second bound \eqref{DF:second_cost_bound} is only preferable to \eqref{DF:first_cost_bound}, if the dimension $d$ is rather small in terms of $L$.

Looking at the three different variants of choices of the frequency set $K=K_\nu$, namely as a  continuous hyperbolic cross (cf. Section~\ref{se:contX}), a $d$-dimensional rectangle (cf. Section~\ref{se:wt_rect}) or a  step hyperbolic cross (cf. Section~\ref{se:stepX}),
%our lattice-based compression algorithm , 
we observe that the minimal precomputation or ``start up'' cost is obtained if we consider the high-dimensional rectangle, i.e., $K = R^{\alpha}_{\nu,d}$. Here, we have a cost of $O(LNd)$, irrespective of the choice of $\nu$.
This is due to the ``Dirichlet kernel trick'': If we sum up 
trigonometric monomials $\omega_\bsk$, all evaluated in the same point $\bsx$, with frequency vectors $\bsk$ from a $d$-dimensional rectangle $R$, then the resulting sum can be reduced to a $d$-fold product of Dirichlet kernels, regardless of the size of $R$, cf. Section~\ref{Subsec:Rectangle_Precomputation}.
Recall that the error bounds we obtained for rectangles are worse than the ones we derived for hyperbolic crosses, which is why the former approach is, in general, less attractive.
Nevertheless, for $\alpha$ large enough (i.e., essentially a bit larger than the dimension $d$), the convergence rate of the errors 
from Corollary~\ref{Cor3.14} for the rectangle setting is higher than the one in \eqref{eq:nets_p2} obtained in \cite{DF2020} with the help of digital $(t,m,d)$-nets. In that regime the rectangle approach has a better error bound as well as a preferable cost bound (in particular, with respect to the dimension dependence) than the approach based on digital nets from \cite{DF2020}.

The lattice-based algorithm using the truncation to a step hyperbolic cross, $K = Q^{\alpha}_{\ell,d}$, exploits that $Q^{\alpha}_{\ell,d}$ is a composition of rectangles. Hence, the ``Dirichlet kernel trick'' from Section~\ref{Subsec:Rectangle_Precomputation} can be used again to get a fast computation of the corresponding compression weights,
see Section~\ref{subsec:precomp_stepX},
which results in a cost of order
\begin{equation*}
    \binom{d - 1 + \ell}{d-1} NLd, \quad \text{ where } \ell := \ceil{\left(\alpha-\half-\sigma \right)\log_2(L)+\frac{d-1}{2}},
\end{equation*}
and $\sigma > 0$ can be chosen sufficiently small. The parameter $\ell$ is chosen according to \RefCol{cor:wtW_stepX}.

The cost of computing the compression weights using \RefAlg{alg:hypX} from Section~\ref{se:precomp_gen_K}, for { $\nu \asymp L^{\alpha-1/2-\sigma(\varepsilon)}$} for the continuous hyperbolic cross $\Knud$, chosen according to \RefCol{cor:wtW_hypX} and \RefCol{cor:wtH_hypX}, along with the estimates from \eqref{eq:cost_general} and \eqref{eq:hypX_prop2}, is of the order
\begin{equation*}
    dN L^{\half - \frac{1}{4\alpha}+\varepsilon} \left( 1 + 2 \zeta(1+2\alpha \varepsilon) \right)^d,
\end{equation*}
where $\varepsilon > 0$ can be chosen arbitrarily small, and the hidden constants are independent of $N$ and $L$, see \RefSec{se:contX}. Essentially, \RefAlg{alg:hypX}  could also be used in combination with the step-hyperbolic cross, and this would cost at most as much as combining it with the continuous hyperbolic cross, since we have shown $\Qmgammad \subseteq \Knugammad$ for $\nu = 2^m$. This approach would improve if there is a faster method at hand to solve the NDFT formulation, see \eqref{eq:ndft}. %which remains an open problem.

As indicated in the analysis for the high-dimensional rectangle, we may get even better cost estimates if we trade for a weaker error bound. 
To illustrate this point, let us discuss an example for $\alpha \geq 2$.
For any $\varepsilon > 0$, there exist $\sigma = \sigma(\varepsilon) > 0$ and $\delta = \delta(\varepsilon)$ such that we obtain for the continuous hyperbolic cross $K = K^\alpha_{\nu,d}$ with
$\nu \asymp L^{\frac{4\alpha^2}{(\alpha+1)(2\alpha-1)}-\sigma}$  from \RefThm{th:hyp_wtH} %then we have again from \RefThm{th:hyp_wtH} 
the error bound
\begin{equation*}
     \err_1(g,\calC) + \err_2(g,\calC) \lesssim \norm{g}_{\Halphad} L^{-\frac{\alpha}{\alpha+1}+  \varepsilon} \log(L)^{\frac{d-1}{2}},
\end{equation*}
and the cost of precomputing the compression weights is of the order
\begin{equation*}
    d N L^{\frac{1}{(\alpha+1)(1-1/2\alpha)}  +  \varepsilon} \left( 1 +   \pi^{-1- 2\alpha \delta} \zeta(1+2\alpha \delta )\right)^d.
\end{equation*}
Note that these results are very close to those obtained in \cite{DF2020} for digital nets.

\subsection*{Acknowledgment}
The authors would like to thank Josef Dick, Michael Feischl, Stefan Kunis, and Dirk Nuyens for valuable discussions.

Part of the work on this paper was done while all three authors attended the Dagstuhl workshop 23351, ``Algorithms and Complexity for Continuous Problems". The authors thank the Leibniz Center for Informatics, Schloss Dagstuhl, and its staff for their support and hospitality. 

The first and third authors would also like to thank the Isaac Newton Institute for Mathematical Sciences (INI), Cambridge, for support and hospitality during the program ``Multivariate approximation, discretization, and sampling recovery'', where the work on this paper was undertaken. This work was supported by an EPSRC grant EP/R014604/1. %Both authors thank the staff of INI for its warm hospitality. 

\section*{Appendix}
\subsection*{Proof of \RefThm{th:wt_rect_wtH}}
To bound the first error term $\err_1(g,\calC)$, we start by estimating the sum appearing in \RefLem{le:trunc_wtH} as
\begin{equation*}
    \begin{aligned}
        \sum_{\bsk \notin \Rnugammad} \frac{1}{r_{\alpha}(\bsgamma, \bsk)} 
        &\leq \sum_{i=1}^d \sum_{\substack{k_i \in \Z\\ r_{\alpha}(\gamma_i, k_i) > \nu}} \frac{1}{r_{\alpha}(\gamma_i, k_i)} \left( 1 + 2 \prod^d_{\substack{j = 1\\ j\neq i}} \sum_{k_j \in \N} \frac{1}{r_{\alpha}(\gamma_j, k_j)} \right) \\
        &= \sum_{i=1}^d \sum_{\substack{k_i \in \N\\ \abs{k_i} > (\gamma_i \nu)^{1/2\alpha}}} \frac{2\gamma_i}{\abs{k_i}^{2\alpha}} \prod_{\substack{j = 1\\ j\neq i}} \left( 1 + 2 \gamma_j \zeta(2\alpha) \right) \\
        &\le \left( 1 + 2 \zeta(2\alpha) \right)^{d-1} \sum_{i=1}^d \sum_{\substack{k_i \in \N\\ \abs{k_i} > (\gamma_i \nu)^{1/2\alpha}}} \frac{2\gamma_i}{\abs{k_i}^{2\alpha}}  \\
        &\lesssim \frac{2}{2\alpha - 1} \left( 1 + 2 \zeta(2\alpha) \right)^{d-1}\sum_{i=1}^d \gamma_i^{1/2\alpha} \nu^{-1+\frac{1}{2\alpha}}\\
        &\lesssim \nu^{-1+\frac{1}{2\alpha}},
    \end{aligned}
\end{equation*}
where in the fourth step we estimated the sum over $k_i$ in a straightforward manner by an integral. 
Now we use estimate~\eqref{eq:error1} and \RefLem{le:trunc_wtH} to conclude
\begin{equation}
    \err_1(g, \calC) \lesssim \norm{g}_{\Halphagammad} \frac{ \mubar_N}{\nu^{\half - \frac{1}{4\alpha} }}.
\end{equation}

To bound the second error term $\err_2(g,\calC)$, we estimate the sum in \RefLem{le:gphi_wtH} as 

\begin{equation*}
    \begin{aligned}
        \sum_{\bsk \in \Rnugammad} r_{\alpha}(\bsgamma, \bsk) 
        &= \sum_{\bsk \in \Rnugammad} \prod_{j \in [d]} r_{\alpha}(\gamma_j, k_j)\\
        &= \prod_{j \in [d]} \sum_{\substack{k_j \in \Z\\ r_{\alpha}(\gamma_j, k_j) \leq \nu}} r_{\alpha}(\gamma_j, k_j) \\
        &= \prod_{j \in [d]} \left( 1 + \frac{2}{\gamma_j} \sum^{\lfloor (\gamma_j \nu)^{1/2\alpha} \rfloor}_{k=1} k^{2\alpha} \right) \\
        &\lesssim \prod_{j \in [d]} \left( 1 + \frac{2\gamma^{1/2\alpha}_j}{2\alpha + 1}\, \nu^{1 + \frac{1}{2\alpha}} \right) \\
        %&\leq 2^d \prod_{j \in [d]} \frac{1}{\gamma_j} \sum_{{k_j \in \N_{0}: k_j \leq \left( \nu \gamma_j \right)^{1/2\alpha} }} k_j^{2\alpha} \\
        %&\lesssim 2^d \prod_{j \in [d]} \frac{1}{\gamma_j} \left( \nu \gamma_j \right)^{(2\alpha + 1)/{2\alpha}}\\ 
        &\lesssim \nu^{\left( 1 + \frac{1}{2\alpha} \right) d},
    \end{aligned}
\end{equation*}
where again in the fourth step we estimated the sum over $k$ by an integral.
Now we use the estimate above and \eqref{eq:lattice_wt_err}, as well as \RefLem{le:gphi_wtH}, to get
\begin{equation}
    \begin{aligned}
        \err_2(g,\calC) &\leq \norm{g \phiRnugamma}_{\Halphaminusgamma} e(\Halphaminusgamma,\calZ) \\
        &\lesssim \norm{g}_{\Halphagammad} \frac{\nu^{\left(\half + \frac{1}{4\alpha} \right)d}}{L^{\alpha-\half-\delta-\tau}} c_{\alpha,\bsgamma,d}\ \zeta_{\delta,d}\ C_{\bsgamma,d}\left(\alpha-\half-\delta,\tau\right)  \mubar_N.
    \end{aligned}
\end{equation}
The statement is therefore proved.

{\small

\printbibliography

}

\end{document}